\documentclass[10pt]{amsart}

\usepackage[top=2cm, bottom=2cm, left=2cm, right=2cm]{geometry}
\usepackage{amsmath, amssymb}
\usepackage{graphicx}
\usepackage{verbatim}
\usepackage[utf8]{inputenc}
\usepackage[colorlinks=true,
            linkcolor=blue,
            citecolor=blue,
            urlcolor=blue]{hyperref}

\usepackage{caption}
\makeatletter
\renewenvironment{abstract}{%
  \ifx\maketitle\relax
    \ClassWarning{\@classname}{Abstract should precede \protect\maketitle}%
  \fi
  \global\setbox\abstractbox=\vbox\bgroup
    \normalfont\small
    \list{}{\labelwidth\z@
      \leftmargin0pt \rightmargin0pt  
      \listparindent\normalparindent \itemindent\z@
      \parsep\z@ \@plus\p@
    }%
    \item[\hskip\labelsep\scshape\abstractname.]%
}{%
  \endlist\egroup
  \ifx\@setabstract\relax \@setabstracta \fi
}
\makeatother

\begin{document}

\title{Density-based structural frameworks for prime numbers, prime Gaps, and Euler Products}
\author{Gregorio Vettori}
\email{gregoriovettori@outlook.com}
\date{January 2026}

\begin{abstract}
We develop a unified density-based framework for primality, coprimality, and prime pairs, and introduce an intrinsic normalized model for prime gaps constrained by the Prime Number Theorem. 
Within this setting, a structural tension between Hardy-Littlewood, Cramér, and PNT predictions emerges, leading to quantitative estimates on the rarity of extreme gaps. 
Additive representations of even integers are reformulated as local density problems, 
yielding non-conjectural upper and lower bounds compatible with Hardy-Littlewood heuristics. 
Finally, the Riemann zeta function is analyzed via truncated Euler products, whose stability and oscillatory structure provide a coherent interpretation of the critical line and prime-based numerical criteria for the localization of non-trivial zeros.
\end{abstract}

\maketitle

\tableofcontents

\section*{Introduction}

In Section~\ref{section 1}, we introduce a unified formalism for primality, coprimality, and prime pairings. Classical arithmetic functions are reinterpreted structurally through normalized density variables, emphasizing collective behavior rather than individual events. This viewpoint provides a natural foundation for the subsequent analysis of prime gaps, additive problems, and Euler products.

Section~\ref{section 2} develops an intrinsic and normalized framework for prime gaps. By introducing a global density function $S_j(n)$ constrained by the Prime Number Theorem, prime gaps are treated as a collective statistical object rather than as independent fluctuations. Within this setting, the Hardy-Littlewood constant $C_2(n)$ emerges as a scale-dependent quantity encoding global correlations, rather than as an ad hoc correction factor. This leads to a clear structural tension between Hardy-Littlewood heuristics, Cramér-type models, and the Prime Number Theorem, and allows for a probabilistic reinterpretation of extreme gaps. Quantitative estimates on the rarity of unusually large gaps follow naturally, without relying on unproven assumptions.

In Section~\ref{section 3}, additive problems are reformulated in terms of local prime densities. Goldbach-type representations are described through a statistical variable $G(n)$, whose behavior is constrained by the same global density principles introduced earlier. Without assuming the Hardy-Littlewood conjecture, we derive non-conjectural upper and lower bounds compatible with it and provide a structural explanation for the extreme improbability of counterexamples. A complementary combinatorial analysis reinforces this picture, showing that potential failures would require highly non-generic and globally correlated exclusions.

Section~\ref{section 4} extends the density-based approach to analytic number theory. 
Within this framework, the critical line $\Re(s)=1/2$ appears as a natural threshold separating stable and unstable regimes of the Euler product. While the analysis is heuristic, it is internally consistent and leads to effective numerical criteria for the localization of non-trivial zeros based solely on prime data. The resulting picture connects naturally with known results on pair correlation and supports the interpretation of non-trivial zeros as global interference phenomena rather than consequences of isolated local behavior.

The overall framework provides a coherent bridge between probabilistic heuristics, combinatorial reasoning, and classical analytic results, and offers a transparent setting in which longstanding conjectures such as Goldbach’s and the Riemann Hypothesis admit a natural structural interpretation.

\section{Explicit formalism for prime numbers}
\label{section 1}

Let $n \geq 2$ be a natural number. In modular arithmetic, an integer $i$ divides $n$ if and only if
$ n \equiv 0 \pmod{i} $.
By definition, prime numbers are those integers $n$ for which this condition holds only for the trivial divisors $i=1$ and $i=n$.

To encode the divisibility structure of integers, we consider all divisions between pairs $(n,i)$ with $2 \leq i \leq n-1$ and organize the corresponding remainders into a triangular array indexed by $(n,i)$. This representation highlights the arithmetic structure underlying divisibility relations.

In this setting, entries corresponding to zero remainder identify non-trivial divisors. For each fixed $i$, these entries form arithmetic progressions corresponding to multiples of $i$. Restricting to $i \leq \lfloor n/2 \rfloor$ captures all non-trivial divisibility relations.

The total number of non-trivial zero-remainder divisions among all integers up to $n$ is therefore given by
\begin{equation}
Z(n) = \sum_{i=2}^{\left\lfloor n/2 \right\rfloor} \left\lfloor \frac{n - i}{i} \right\rfloor = 1 + \sum_{i=3}^{\left\lfloor n/2 \right\rfloor} \left\lfloor \frac{n}{i} \right\rfloor \\
= 1 - n - \left\lfloor \frac{n}{2} \right\rfloor + \sum_{i=1}^{\left\lfloor n/2 \right\rfloor} \left\lfloor \frac{n}{i} \right\rfloor,
\label{Z}
\end{equation}
where the last expression is valid for all $n \geq 2$.

The summatory function $Z(n)$ does not introduce new arithmetic information,
but rather reorganizes classical divisor data into a form that is well-suited
for normalization and indicator constructions.

The number of non-trivial positive divisors of an integer \(n\) can be expressed in terms of the summatory function \(Z(n)\) as
\begin{equation}
P(n)=Z(n)-Z(n-1)
= \sum_{i=2}^{\lfloor n/2 \rfloor} \left\lfloor \frac{n-i}{i} \right\rfloor
     - \sum_{i=2}^{\lfloor (n-1)/2 \rfloor} \left\lfloor \frac{n-1-i}{i} \right\rfloor,
\label{P(n)}
\end{equation}
where \(P(n)\) counts the number of integers \(i\in(1,n)\) such that \(i\mid n\).
By construction, \(Z(n)=\sum_{k=4}^{n}P(k)\), with \(Z(1)=Z(2)=Z(3)=0\).

Since an integer \(n\) is prime if and only if it admits no non-trivial divisors, the condition \(P(n)=0\) provides a necessary and sufficient criterion for primality, consistent with classical divisor-based characterizations.

Exploiting the properties of the floor function and the fact that all primes greater than \(2\) are odd, for odd \(n\) the expression for \(P(n)\) simplifies to
\begin{equation}
P(n) = \sum_{i=2}^{(n-1)/2}
\left(
\left\lfloor \frac{n}{i}\right\rfloor - \left\lfloor \frac{n-1}{i} \right\rfloor \right).
\end{equation}

$\lfloor (n-2)/2 \rfloor$ reflects the maximal number
of admissible non-trivial divisors for an integer $n$.
We define
\begin{equation}
\omega(n) = 1 - P(n)\,\left\lfloor \frac{n-2}{2} \right\rfloor^{-1},
\end{equation}
as the fraction of divisions yielding a non-zero remainder with respect to the maximal theoretical case.
The normalization by $\lfloor (n-2)/2 \rfloor$ ensures that $\omega(n)$ takes values in $[0,1]$.
Within this scale, prime numbers correspond to the extremal case $\omega(n)=1$,
while highly composite integers occupy lower values.
In particular, \(\omega(n)=0\) occurs only for \(n=4\) and \(n=6\).

Within this framework, an explicit prime indicator function can be defined for \(n\geq 4\) as
\begin{equation}
I(n)=\left\lfloor \omega(n)\right\rfloor =
\begin{cases}
1, & \text{if } n \text{ is prime},\\
0, & \text{otherwise}.
\end{cases}
\label{I}
\end{equation}
By construction, \(I(n)=1\) if and only if no non-trivial congruence
\(n\equiv 0 \pmod{i}\) admits a solution for \(i\in[2,\lfloor n/2\rfloor]\).
This provides an explicit arithmetic characterization of primality within the proposed formalism.

The number of primes less than or equal to \(n\) is then given by
\begin{equation}
K(n)=2+\sum_{j=4}^{n} I(j),
\label{K}
\end{equation}
where the initial offset accounts for the primes \(2\) and \(3\), and
\(I(n)=K(n)-K(n-1)\).
This construction parallels the classical definition of the prime counting function as a cumulative sum of a prime indicator \cite{HW}.

The functions \(P(n)\) and \(I(n)\) are directly related to the classical divisor function \(\tau(n)\) by
\begin{equation}
P(n) = \tau(n)-2, \qquad
I(n) = \left\lfloor 1 - \big(\tau(n)-2\big)\left\lfloor \frac{n-2}{2}\right\rfloor^{-1} \right\rfloor,
\end{equation}
where \(\tau(n)\) denotes the total number of positive divisors of \(n\).
For an integer factorization \(n=p_1^{q_1}p_2^{q_2}\cdots p_k^{q_k}\), one has
\(\tau(n)=(q_1+1)(q_2+1)\cdots(q_k+1)\).

From a computational perspective, the present formulation is not intended to compete with classical sieve-based methods, such as the Sieve of Eratosthenes or the Sieve of Atkin, for the efficient enumeration of prime numbers.
Its purpose is instead to provide an exact and explicit arithmetic representation of primality directly derived from divisor properties.

Similarly, when interpreted as a criterion for primality, the construction is not optimized for large-scale numerical testing and is therefore not comparable, in terms of computational complexity, to established probabilistic or deterministic methods such as Fermat-type tests, Wilson's theorem, or the AKS algorithm.
The interest of the present approach lies in its structural transparency and in the fact that it yields a closed-form characterization without relying on probabilistic assumptions or auxiliary algebraic constructions.

In the following sections, this formalism is employed as a unified framework to investigate structural properties of prime numbers and their gaps, and to examine its implications for several classical problems in analytic number theory, including conjectures of Hardy–Littlewood type, Goldbach-type representations, and properties of the Riemann zeta function.

\subsection{Prime numbers density and distribution}

The prime density in the interval $[1,n]$ is defined as $D(n)=K(n)/n$, and the corresponding average spacing between consecutive primes is $L(n)=1/D(n)$.
A local version of these quantities can be introduced by considering a symmetric window of width $2x$ around $n$:
\begin{equation}
D_x(n) = \frac{K(n+x)-K(n-x)}{2x}, \qquad
L_x(n) = \frac{1}{D_x(n)}.
\label{D-L}
\end{equation}

According to the Prime Number Theorem, the prime counting function satisfies
\begin{equation}
K(n)\sim \mathrm{Li}(n)=\int_2^n\frac{du}{\ln u}\sim \frac{n}{\ln n},
\end{equation}
where $\mathrm{Li}(n)$ provides a well-known asymptotic approximation \cite{MV}.

We introduce the following heuristic refinement of $n/\ln n$,
\begin{equation}
A(n)=\frac{n}{\ln n}\left(1+\frac{1.08}{(\ln n)^{1.01}}\right),
\label{new_K(n)_try}
\end{equation}
which is fully consistent with the Prime Number Theorem and its classical refinements.
In particular, Table~\ref{table:5} shows that $A(n)$ provides a significantly closer approximation to $\mathrm{Li}(n)$ than $n/\ln n$ over the explored range.

\begin{table}[h!]
\centering
\begin{tabular}{l c c} 
 \hline
 $n$          & $n / (\text{Li}(n) \ln n) - 1$ & $A(n)/\text{Li}(n) - 1$  \\ [1ex]
 \hline
 $10^{10^1}$         & $-4.56 \cdot 10^{-2}$  & $-2.24 \cdot 10^{-3}$   \\
 $10^{10^2}$         & $-4.36 \cdot 10^{-3}$  & $-6.06 \cdot 10^{-5}$   \\
 $10^{10^3}$         & $-4.35 \cdot 10^{-4}$  & $-5.76 \cdot 10^{-7}$   \\
 $10^{10^4}$         & $-4.34 \cdot 10^{-5}$  & $-1.01 \cdot 10^{-6}$    \\
 $10^{10^5}$         & $-4.34 \cdot 10^{-6}$  & $-1.97 \cdot 10^{-7}$    \\
 $10^{10^6}$         & $-4.34 \cdot 10^{-7}$  & $-2.91 \cdot 10^{-8}$   \\
 \hline
\end{tabular}
\caption{Relative errors of the approximations $n/\ln n$ and $A(n)$ with respect to the logarithmic integral $\mathrm{Li}(n)$.}
\label{table:5}
\end{table}

\subsection{Average and expected value of primes}

In connection with the distribution of primes, it is natural to consider the average value of all prime numbers not exceeding a given integer \(n\).
Within the present framework, this quantity can be defined as
\begin{equation}
\bar p(n) = \frac{1}{K(n)} \sum_{i=1}^{n} i\, I(i),
\label{p_n}
\end{equation}
where \(I(i)\) denotes the prime indicator introduced in Section \ref{section 1}.

The ratio \(\bar p(n)/n\) is strictly smaller than \(1/2\) for moderately large values of \(n\), and numerical evidence suggests the presence of a shallow minimum at small \(n\).
If prime numbers were uniformly distributed in \([1,n]\), one would instead obtain the trivial estimate \(\bar p(n)=n/2\).

\begin{figure}[h!]
\begin{center}
  \includegraphics[width=0.5\columnwidth]{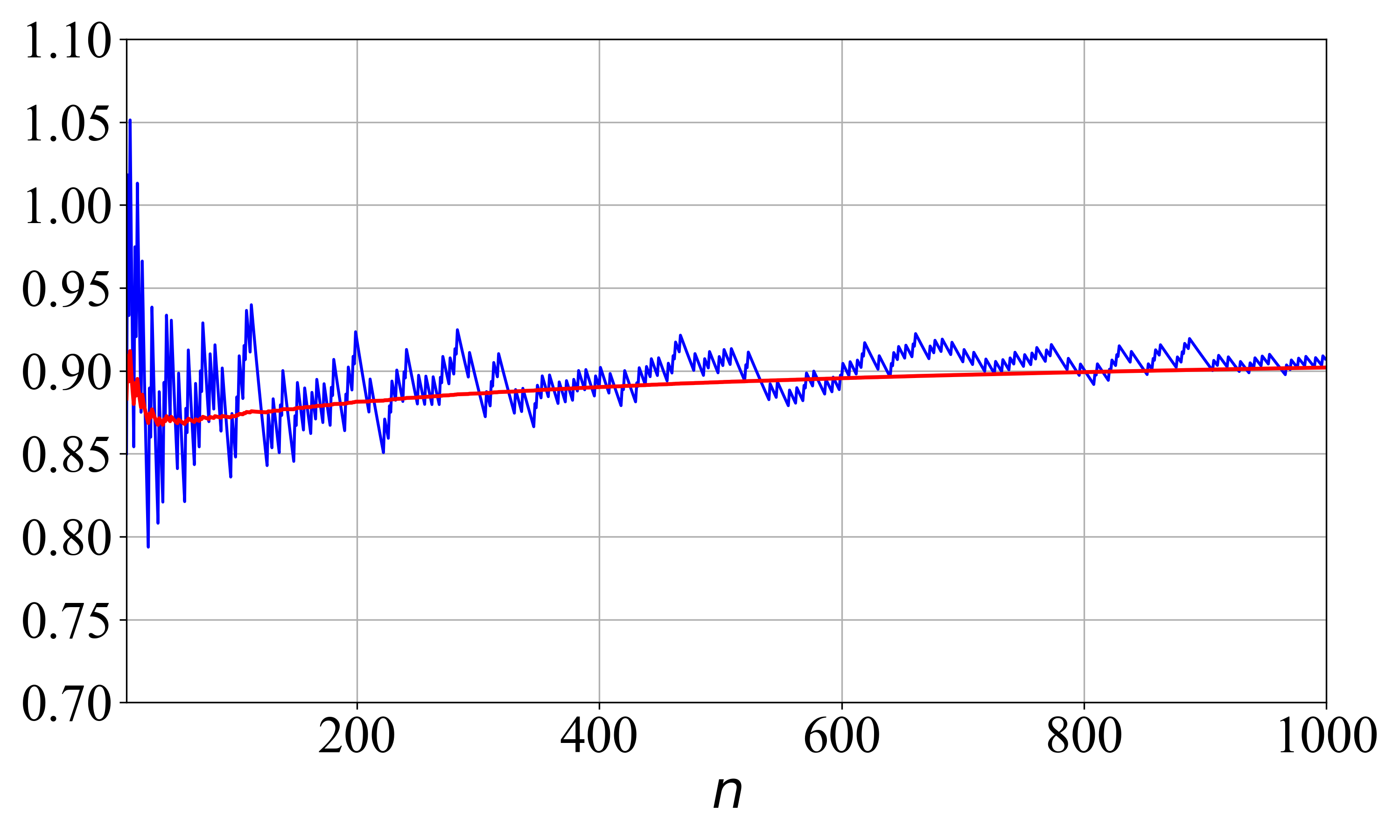}
  \caption{$2 \ \bar p(n)/n$ \eqref{p_n} (blue) and its heuristic approximation \eqref{p_avg} (red), for $n\leq 1000$.}
  \label{p average}
\end{center}
\end{figure}

The deviation from this uniform behavior can be heuristically understood in terms of the decreasing prime density.
Indeed, using the Prime Number Theorem, the ratio between the average densities in the intervals \([1,n]\) and \([1,n/2]\) satisfies
\begin{equation}
\frac{D(n)}{D(n/2)} \sim \frac{\ln(n/2)}{\ln n} < 1,
\end{equation}
which implies that primes are, on average, more concentrated in the lower half of the interval.
This leads to the heuristic approximation
\begin{equation}
\bar p(n) \sim \frac{n}{2} \frac{\ln(n/2)}{\ln n} + K(n)\frac{\ln n}{n},
\label{p_avg}
\end{equation}
which tends to \(n/2\) as \(n\to\infty\) (Figure \ref{p average}), consistently with the asymptotic uniformity predicted by the Prime Number Theorem \cite{MV}.

\subsection{Coprimality and structural relations}

Two integers are said to be \emph{coprime} if their only common positive divisor is \(1\).
Within the present framework, the number of common divisors $ > 1$ for two integers \(n\) and \(x\) can be written as
\begin{equation}
C(n,x) =
\sum_{i=2}^{\left\lfloor \min(n,x)/2 \right\rfloor}
\left(
\left\lfloor \frac{n}{i} \right\rfloor -
\left\lfloor \frac{n-1}{i} \right\rfloor
\right)
\left(
\left\lfloor \frac{x}{i} \right\rfloor -
\left\lfloor \frac{x-1}{i} \right\rfloor
\right).
\end{equation}

By construction, \(C(n,x)=0\) if and only if \(n\) and \(x\) share no common divisors greater than \(1\), and hence are coprime.
In particular, an integer \(p\ge2\) is prime if and only if \(C(n,p)=0\) for all integers \(n\) such that \(1<n<p\).

The maximum theoretical number of common divisors greater than \(1\) is bounded above by \(\lfloor x/2 \rfloor\).
This motivates the introduction of the coprime indicator function
\begin{equation}
I_C(n,x) =
\left\lfloor
1 - C(n,x)\,\left\lfloor \frac{x}{2} \right\rfloor^{-1}
\right\rfloor,
\label{I_C}
\end{equation}
where \(I_C(n,x)=1\) if \(n\) and \(x\) are coprime and \(I_C(n,x)=0\) otherwise.

The total number of integers less than or equal to \(n\) that are coprime with \(x\) is then given by
\begin{equation}
K_C(n,x) = \sum_{i=2}^{n} I_C(i,x).
\label{K_C(n,n)_eq}
\end{equation}
In the special case \(n=x\), one has
\[
K_C(n,n)=\varphi(n) = n \prod_{p|n} \bigg( 1 - \frac{1}{p} \bigg),
\]
where \(\varphi(n)\) denotes Euler's totient function and $p_i$ are all the primes that are part of the factorization of $n$.

It is well known that the average order of the totient function satisfies
\[
\frac{\varphi(n)}{n} \to \frac{6}{\pi^2}
\quad \text{on average},
\]
which explains the typical density of integers coprime with \(n\).

The maximal value of $\varphi(n)$ is attained when $n=p$ is prime, while the
minimum value occurs for integers having the largest possible number of small
prime factors. In particular, fixing a prime $w \in \mathbb{P}$, for all integers
$n < n_{\max} = \prod_{p \leq w} p$ one has
\begin{equation}
\varphi(n)_{\min} = n \prod_{p \leq w} (1 - p^{-1}), 
\qquad 
\varphi(n)_{\max} = n (1 - n^{-1}),
\label{K_C_bounds}
\end{equation}
corresponding respectively to minimal and maximal coprimality constraints.

Beyond these extremal bounds, the values of $K_C(n,n)$ exhibit a pronounced
statistical clustering. In particular, $K_C(n,n)$ frequently concentrates near
the values $n/3$, $n/2$, and $2n/3$. This behavior reflects the dominant influence
of the smallest prime divisors, especially $2$ and $3$, on coprimality densities.

\begin{figure}[h!]
\begin{center}
  \includegraphics[width=0.5\columnwidth]{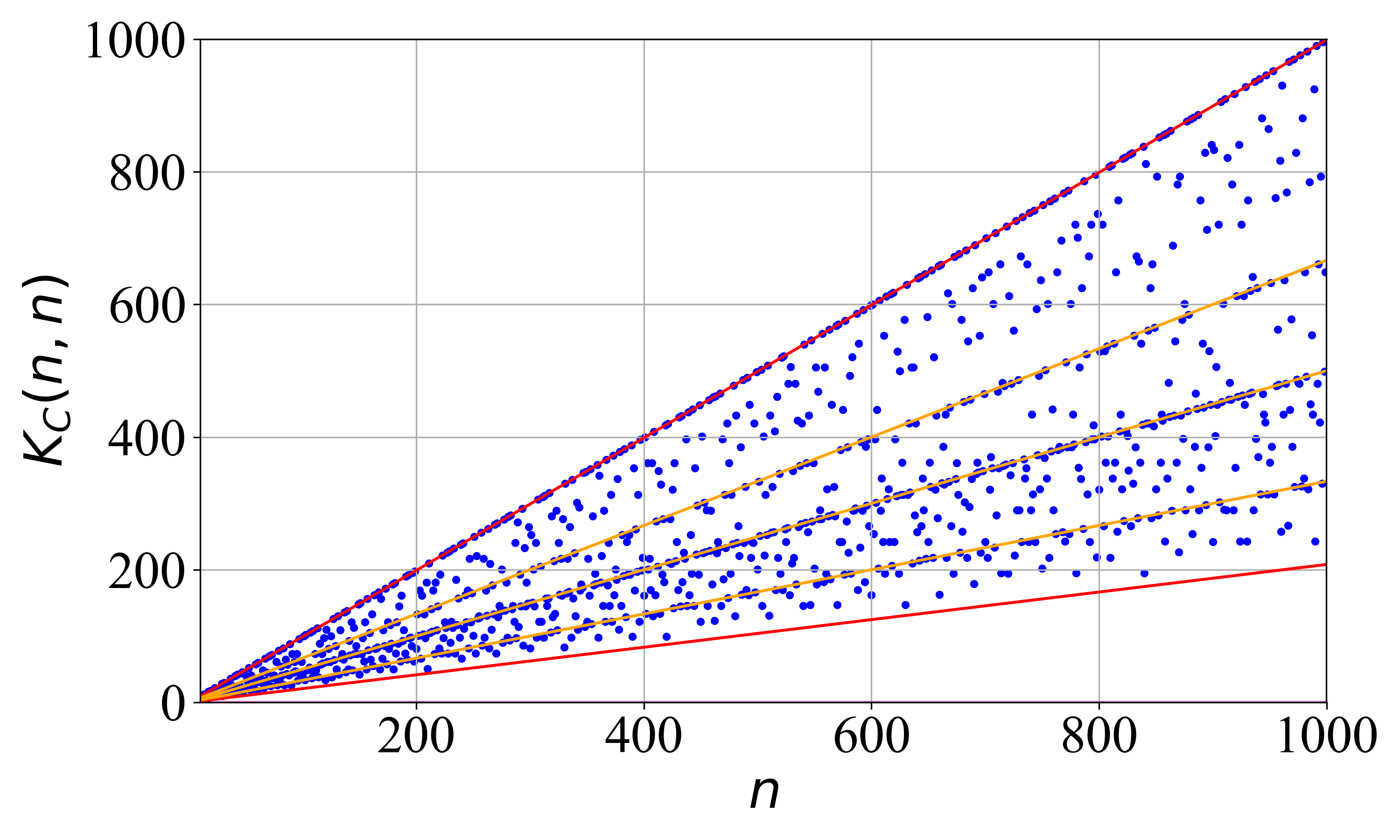}
  \caption{$K_C(n,n)$ \eqref{K_C(n,n)_eq} (blue), $\varphi(n)_{min} \sim 0.208 \ n$ and $\varphi(n)_{max}$ \eqref{K_C_bounds} (red), for $n \leq 1000$ and $w = 11$. As reference, $n/3$, $n/2$ and $2n/3$ (orange).}
  \label{co-primes}
\end{center}
\end{figure}

\begin{figure}[h!]
\begin{center}
  \includegraphics[width=0.5\columnwidth]{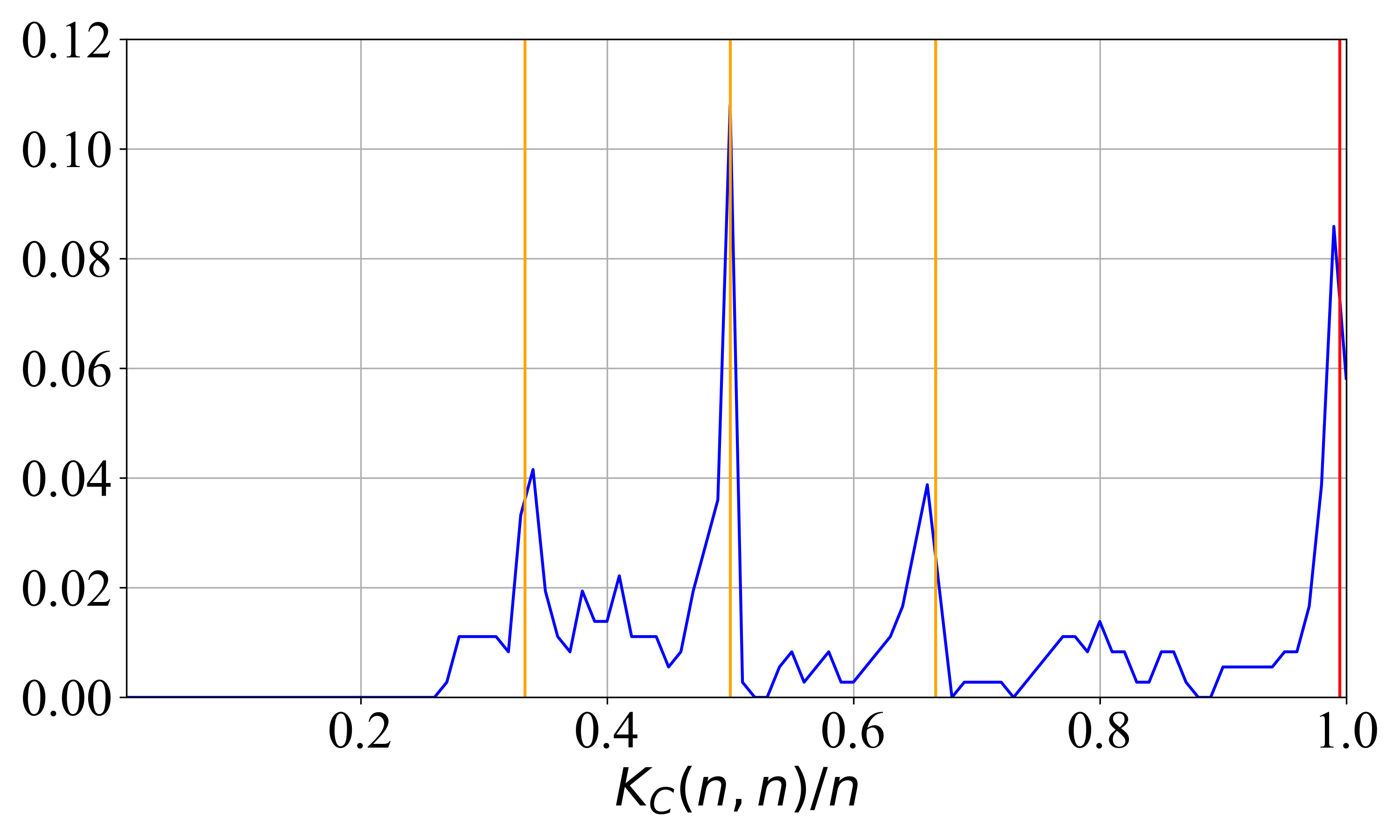}
  \caption{Distribution of $K_C(n,n)/n$ values (blue), for $n \leq 200$. As reference, $1-n^{-1}$ (red), $1/3$, $1/2$ and $2/3$ (orange).}
  \label{co-primes-dist}
\end{center}
\end{figure}

This phenomenon is consistent with the arithmetic constraint that all primes
greater than $3$ lie in the residue classes $6k \pm 1$, and highlights the role
of low-modulus structures in shaping global coprimality statistics. 
The effect
is illustrated in Figure~\ref{co-primes}, which displays $K_C(n,n)$ for $n \leq
1000$, and in Figure~\ref{co-primes-dist}, showing the empirical distribution of
$K_C(n,n)/n$ for $n \leq 200$.

\subsection{Twin primes and prime pairs}

A pair of primes \((p,p+2)\) is called a \emph{twin prime pair}.
Within the present framework, an explicit indicator function for twin primes can be defined for \(n\ge4\) as
\begin{equation}
I_2(n) =
\left\lfloor
\frac{1}{2}
\left(
I(n)\big(I(n-2)+I(n+2)\big) + 1
\right)
\right\rfloor.
\label{I2}
\end{equation}

The symmetric and centered form of \(I_2(n)\) is intentional.
It allows one to identify individual primes belonging to a twin prime pair, independently of whether the companion prime occurs at distance \(+2\) or \(-2\).
In this way, the indicator is associated with single prime elements rather than with unordered prime pairs.

More generally, for any even integer \(x\ge2\), we define an indicator function for prime pairs separated by a distance \(x\) as
\begin{equation}
I_x(n) =
\left\lfloor
\frac{1}{2}
\left(
I(n)\big(I(n-x)+I(n+x)\big) + 1
\right)
\right\rfloor.
\label{I_X}
\end{equation}

The total number of prime pairs with separation \(x\) and with elements not exceeding \(n\) is then given by
\begin{equation}
K_x(n) = \gamma_x + \sum_{j=4}^{n} I_x(j),
\end{equation}
where the correction term \(\gamma_x = I(3)\,I(3+x)\) accounts for boundary effects arising from the lower limit of the summation.
This correction is negligible in the asymptotic regime.

The indicator functions $I_x(n)$ provide an explicit arithmetic realization
of the prime pair counting functions appearing in the Hardy-Littlewood conjectures
\cite{HW1923}, without assuming their conjectural asymptotic formulas.

By construction, \(I_x(n)\) takes values only in \(\{0,1\}\), even in the presence of prime triplets or higher-order clusters.
Moreover, \(K_x(n)\) counts all prime pairs separated by \(x\), independently of whether they are consecutive.
For \(x\ge6\), this leads to an overcounting of consecutive prime gaps of size \(x\), since intermediate primes may exist.
In contrast, for \(x=2\), such overcounting is impossible, and \(K_2(n)\) yields the exact counting function for twin prime pairs.

The formulation of prime pairs within the same indicator framework
used for primality and coprimality highlights the structural unity
of the present approach and prepares the density-based analysis
of additive representations developed in Section \ref{section 3}.

\newpage

\section{An intrinsic framework for prime gaps}
\label{section 2}

The Hardy-Littlewood conjecture predicts that the number of twin prime pairs not exceeding \(n\) satisfies the asymptotic relation
\begin{equation}
K_2(n) \equiv 1 + \sum_{j=4}^{n} I_2(j) \sim 2 C_2 \int_2^n \frac{dt}{(\ln t)^2}
          \sim \frac{2 C_2\, n}{(\ln n)^2},
\end{equation}
where \(C_2\) denotes the twin prime constant,
\begin{equation}
C_2 = \prod_{p \geq 3} \left( 1 - \frac{1}{(p-1)^2} \right) \approx 0.660161.
\end{equation}

Motivated by this conjecture, we consider the problem in a more general setting by studying the distribution of prime gaps of fixed even size \(j\) among primes not exceeding \(n\), with the aim of revealing structural constraints in the spacing of prime numbers.

A classical heuristic model introduced by Cramér suggests that the maximal size of prime gaps grows on the order of \((\ln n)^2\) in the asymptotic regime \cite{Cramer1936}.
Accordingly, we restrict attention to even gap sizes \(j\) satisfying
\[
2 \leq j \leq (\ln n)^2,
\]
a range which includes twin primes and excludes exceptionally large gaps.

Empirically, counting the prime gaps equal to \(j\) among all primes not exceeding \(n\) yields approximately \(n/\ln n\) observations, corresponding to the number of primes up to \(n\).
Normalizing the frequency of each gap size by \(n\), we define \(P_{\mathrm{emp}}(j,n)\) as the empirical density of prime gaps of size \(j\).

Under a naive random model in which primality occurs independently with probability \(1/\ln n\), the expected density of pairs of integers at distance \(j\) both being prime is of order \(1/(\ln n)^2\).
Observed prime gap distributions deviate significantly from this baseline and exhibit a strong positive skewness.
While the mean prime gap is of order \(\ln n\), the majority of observed gaps are substantially smaller, with small gap sizes occurring with much higher frequency than larger ones.

This asymmetry and the presence of a long tail toward large gaps will be analyzed in detail in the following subsections.
In particular, the present approach aligns with the modern view of primes as governed by an interplay between structure and randomness, as articulated for instance by Tao \cite{Tao2015}, while emphasizing the role of global normalization constraints.

\subsection{Intrinsic model for prime gaps}

We introduce a heuristic model for the empirical distribution of prime gaps, designed to capture the main statistical features observed in numerical data, including asymmetry, characteristic scale, and arithmetic modulation.

Specifically, we model the probability density of prime gaps \(j\) by a log-normal distribution,
\begin{equation}
P_j(n) =
\frac{C(n)}{\sigma \sqrt{2\pi}\, e^x}
\exp\!\left(-\frac{(x-\mu)^2}{2\sigma^2}\right),
\end{equation}
where \(x=\ln j\), with \(j\in(2,(\ln n)^2)\),
\(\mu=\ln(\ln n)\), and \(\sigma=\mu/\alpha(j,n)\).

The choice of a log-normal profile is not meant to be canonical.
Rather, it reflects two minimal requirements simultaneously observed
in empirical prime gap data: strong right-skewness with a long tail and a characteristic scale of order $\ln n$.
Any smooth model satisfying these constraints would lead to qualitatively
similar conclusions in the subsequent analysis.
Similar log-normal heuristics for prime gaps have appeared in refinements of Cramér-type models \cite{Granville1995}.

The distribution is discrete, as admissible prime gaps are even integers, and is normalized through the factor \(C(n)\) so that
\(\sum_j P_j(n) \sim 1/\ln n\).

Empirically, gaps \(j\) divisible by \(6\) occur with enhanced frequency.
This effect reflects the arithmetic constraint that all primes greater than \(3\) lie in the residue classes \(6k\pm1\).
To account for this structure, we introduce a multiplicative modulation with period \(6\),
\begin{equation}
S_j(n) \sim \frac{P_j(n)}{2}\left[3+\cos\!\left(\frac{\pi j}{3}\right)\right],
\label{S_j(n)}
\end{equation}
yielding \(S_j(n)\in(P_j(n),2P_j(n))\).

Although \(S_j(n)\) is not normalized, it provides a significantly improved description of the empirical gap frequencies.
To quantify this improvement, we define the mean relative deviation
\begin{equation}
E(n)= (\ln n)^{-2}
\sum_{j=2}^{\lfloor(\ln n)^2\rfloor}
\left|\frac{S_j(n)}{P_{\mathrm{emp}}(j,n)}-1\right|.
\end{equation}
Numerical experiments up to $n=10^7$ indicate that $E(n)$ decreases with $n$.

Motivated by these observations, we consider a generalized form
\begin{equation}
S_j(n) \sim
\frac{\beta(n)}{\mu e^{\mu(\rho+1)}}
\frac{3+\cos(\pi j/3)}
{2\exp\!\left(\tfrac{1}{2}(\rho-1)^2\alpha(j,n)^2\right)},
\end{equation}
where \(\rho=(\ln j)/\mu\) and \(\alpha(j,n)\) is allowed to vary with \(j\).
The variable \(\rho\), normalized to the interval \((2/\ln n,2)\), provides a more natural scaling parameter for prime gaps than the canonic merit \(j/\ln n\).

\subsection{PNT global density constraints}

We impose six asymptotic constraints on the function \(S_j(n)\), reflecting known properties of prime gaps, normalization requirements, and established conjectural results.
The constraints are not assumed to hold simultaneously
in the true distribution of prime gaps.
Instead, they represent structural requirements imposed on $S_j(n)$
to test the internal compatibility of classical conjectures
within a single normalized framework.
In particular, the constraints provide information on the asymptotic behavior of the functions \(\beta(n)\) and \(\alpha(j,n)\) for \(n\gg1\).

\begin{subequations}
\begin{align}
& S_{\ln n}(n) \sim (\ln n)^{-2} \\
& S_2(n) \gtrsim 2C_2\,(\ln n)^{-2} \\
& S_{(\ln n)^2}(n) \lesssim n^{-1}(\ln n)^{-2} \\
& S_4(n) \sim S_2(n) \\
& n \sum_{j\le (\ln n)^2} j\, S_j(n) \sim n \\
& n \sum_{j\le (\ln n)^2} S_j(n) \sim \mathrm{Li}(n)
\end{align}
\label{constraints}
\end{subequations}

For \(\rho\sim1\), the function \(S_j(n)\) is effectively independent of \(\alpha(j,n)\).
The first constraint fixes the overall normalization, yielding \(\beta(n)\sim 2\mu/3\), and
\begin{equation}
S_j(n) \sim N_j
\exp\!\left[-\mu(\rho+1)-\tfrac12(\rho-1)^2\alpha(j,n)^2\right],
\end{equation}
where \(N_j = 1+\cos(\pi j/3)/3 \in (0.83,1.33)\) and \(\langle N_j\rangle=1\).

For fixed gap \(j=2\) and \(n\to\infty\), we assume that the density of prime gaps is governed by the Hardy–Littlewood conjecture.
Imposing the second constraint yields
\begin{equation}
\frac12\,\alpha^2(2,n)
\sim
\frac{\mu}{1-\rho}
-
\frac{\ln(12C_2/5)}{(1-\rho)^2}
\;\to\;
\mu,
\end{equation}
where \(C_2\simeq0.660161\) is the twin prime constant.

For \(\rho\sim2\), prime gaps become extremely rare.
The third constraint implies
\begin{equation}
\frac12\,\alpha^2((\ln n)^2,n)
\sim
\frac{\ln n-\mu(\rho+1)}{(\rho-1)^2}
\;\to\;
\ln n,
\end{equation}
effectively determining the decay rate of the tail of \(S_j(n)\).

The fifth and sixth constraints ensure the correct total number of prime gaps, $ \sim \text{Li}(n)$, and an average gap size of order \(\ln n\).
The fourth constraint, $S_2(n)\sim S_4(n)$,
is empirical and reflects the observed similarity between the frequencies of gaps of size \(2\) and \(4\).
It yields
\begin{equation}
S_2(n)\sim S_4(n)\sim
\frac12\left(
\frac{\mathrm{Li}(n)}{n}
-
\sum_{j=6}^{\lfloor(\ln n)^2\rfloor} S_j(n)
\right).
\end{equation}

We generalize the previous analysis by introducing two linear dependencies in the normalized variable
\(\rho = \ln j / \mu\), corresponding to the regions below and above the average prime gap scale.

Rather than prescribing a single functional form for the effective parameter \(\alpha(j,n)\) over the entire range of \(\rho\),
I interpolate directly the resulting density \(S_j(n)\) in the region \(\rho \in (\rho_2,1)\),
where the behavior is controlled by small and intermediate gaps,
while retaining an explicit parametrization in terms of \(\alpha(j,n)\) for \(\rho \in (1,2)\),
where the tail of the distribution becomes dominant.

For \(\rho \in (\rho_2,1)\), we assume the linear relation
\begin{equation}
S_j(n)\, e^{2\mu} N_j^{-1}
\sim
\frac{12}{5} C_2 + \Bigl(1 - \frac{12}{5} C_2\Bigr)(\rho - \rho_2),
\label{alpha_1}
\end{equation}
which interpolates between the Hardy-Littlewood prediction for twin primes at \(\rho=\rho_2\)
and the average gap scale at \(\rho \simeq 1\).

For \(\rho \in (1,2)\), we instead prescribe a linear growth for the effective width parameter,
\begin{equation}
\frac{1}{2}\,\alpha^2(j,n)
\sim
\mu + (\ln n - 2\mu)(\rho - 1),
\label{alpha_2}
\end{equation}
ensuring a progressive suppression of large prime gaps as \(j\) approaches the maximal scale \((\ln n)^2\).

With these assumptions, the density satisfies the correct asymptotic behaviors:
\(S_j(n) \sim 2C_2\,e^{-2\mu}\) at \(\rho=\rho_2\),
\(S_j(n) \sim e^{-2\mu}\) near \(\rho \simeq 1\),
and \(S_j(n) \sim e^{-2\mu}\,n^{-1}\) as \(\rho \to 2\).

Accordingly, the density \(S_j(n)\) admits two asymptotically consistent representations:
\begin{subequations}
\begin{align}
&S_j(n) \sim
N_j
\Bigl[
\frac{12}{5} C_2 + \Bigl(1 - \frac{12}{5} C_2\Bigr)(\rho - \rho_2)
\Bigr]
e^{-2\mu},
\\[1ex]
&S_j(n) \sim
N_j \ e^{
-\mu(\rho+1)
- (\rho-1)^2\bigl[\mu + (\ln n - 2\mu)(\rho - 1)\bigr]},
\end{align}
\label{final_Sj(n)}
\end{subequations}
with the former valid for \(\rho \in (\rho_2,1)\) and the latter for \(\rho \in (1,2)\).

The quantitative validity of the proposed framework will be assessed through direct comparison with empirical prime gap data.
In particular, the scaling behavior of \(S_j(n)\) with respect to both \(j\) and \(n\),
as well as its accuracy in reproducing the observed distribution of small and intermediate gaps,
will be examined in the following subsections.

Special attention will be devoted to the dependence of the effective parameter \(\alpha(j,n)\)
on the normalized variable \(\rho=\ln j/\mu\),
and to the stability of the resulting fit as \(n\) increases.

The model is heuristic in nature and is not intended to capture fine-scale arithmetic correlations
or to provide exact predictions for individual gap occurrences.
Rather, it is designed to describe the global statistical structure of prime gaps in an averaged sense.

\subsection{Structural tension between HL, Cramér, and PNT}

The sixth constraint in Equation \eqref{constraints} encodes the Prime Number Theorem,
requiring that the total number of prime gaps up to scale \((\ln n)^2\)
be asymptotically equivalent to \(\mathrm{Li}(n)\).
To investigate the internal consistency of the proposed model,
it is therefore natural to introduce a set of aggregate quantities
associated with the density \(S_j(n)\).

\begin{itemize}

\begin{item}
    the expected number of prime gaps equal to $j$ in the range $(n_1,n_2)$ is
\begin{equation}
\int_{n_1}^{n_2} n S_j(n) \ dn \sim n_2 S_j(n_2) - n_1 S_j(n_1)
\end{equation}
\end{item}

\begin{item}
    the expected number of prime gaps equal to $j \in (j_1, j_2)$ is $\sum_{j_1}^{j_2} nS_j(n)$
\end{item}

\begin{item}
    $\varepsilon(n)$ is the error of $n S_j(n)$ with respect to $\text{Li}(n)$:
\begin{equation}
\varepsilon(n) = \frac{\sum_2^{\lfloor j_{max} \rfloor} nS_j(n)}{\text{Li}(n)} - 1
\end{equation}
\end{item}

\begin{item}
The quantity \(k_1(n)\) measures the proportion of prime gaps
whose size lies below the average scale \(\ln n\),
thus providing a natural indicator of the relative weight
of small gaps within the overall distribution:
\begin{equation}
k_1(n) \sim \frac{\sum_{j=2}^{\lfloor \ln n \rfloor} S_j(n)}{\sum_{j=2}^{\lfloor j_{max} \rfloor} S_j(n)}
\label{k1}
\end{equation}
\end{item}

\begin{item}
the average density of prime gaps with $\rho < 1$ is
\begin{equation}
\langle S_j(n) \rangle_1 \sim 2 \ k_1(n) \ (\ln n)^{-2}
\label{Sjn1}
\end{equation}
\end{item}

\begin{item}
the average density of prime gaps with $\rho > 1$ is
\begin{equation}
\langle S_j(n) \rangle_2 \sim \frac{2 \ k_2(n) \ (\ln n)^{-1}}{(\ln n)^{\rho_{max}} - \ln n} \sim \frac{2 \ k_2(n)}{(\ln n)^{\rho_{max}+1}}
\label{Sjn2}
\end{equation}
where $k_2(n) = 1 - k_1(n)$ and $\rho_{max}$ is the upper limit of $\rho$ ($\rho_{max} = 2$ according to the Cramér conjecture).
\end{item}

\end{itemize}

Equations \eqref{final_Sj(n)} and \eqref{k1} lead to $\langle S_j(n) \rangle_1 \sim (C_2(n)+1/2) \ (\ln n)^{-2}$.
Using Equation \eqref{Sjn1}, we obtain an intrinsic relation between the fraction of small gaps \(k_1(n)\)
and an effective twin-prime constant \(C_2(n)\),
which emerges here as a derived quantity, constrained by global normalization rather than imposed a priori.

\begin{equation}
k_1(n) \sim \frac{1}{2}C_2(n) + \frac{1}{4}, \qquad
C_2(n) \sim 2k_1(n) - \frac{1}{2}.
\label{k1-C2}
\end{equation}

Equation \eqref{k1-C2} should be interpreted as a consistency relation
within the present framework.
It does not assert that the twin-prime constant determines the global
gap distribution in reality, but rather that any density model satisfying
the imposed normalization constraints must exhibit such a correlation.

Since \(2C_2(n)\) must be greater than \(1\), it follows that \(k_1(n)\) is always greater than \(50\%\).
Conversely, the constraint \(k_1(n)<1\) implies \(C_2(n)<1.5\).
These bounds follow from general normalization and positivity constraints
and do not rely on the validity of either the Hardy-Littlewood or Cramér conjectures.
In particular, assuming the Hardy-Littlewood prediction
\(C_2(n)\to C_2\simeq 0.66\) would imply \(k_1(n)\to 0.58\) as \(n\to\infty\),
a value consistent with the predominance of gaps below the average scale.

We generalize Equation \eqref{final_Sj(n)} by considering $C_2 = C_2(n)$ (no Hardy-Littlewood conjecture), $\beta(n) = \mu \ r(n)$ (variable normalization), and $\rho_{max} = \rho_{max}(n)$ (no Cramér conjecture), obtaining:

\begin{subequations}
\begin{align}
&S_j(n) \sim y(j,n) \ \big[ z(n) + (1 - z(n))(\rho - \rho_2) \big] \ e^{-2\mu} \\
&S_j(n) \sim y(j,n) \ e^{-\mu (\rho+1) - (\rho-1)^2 \alpha^2(j,n)/2}
\end{align}
\label{final_Sj(n)_2}
\end{subequations}
where $y(j,n) = 3 r(n) N_j /2$, $z(n) = 8C_2(n) / 5 r(n)$, and
\begin{equation}
\frac{1}{2}\alpha^2(j,n) = \frac{\mu (\rho_{max}-1) + (\ln n - \mu\rho_{max})(\rho-1)}{(\rho_{max}-1)^{3}}
\end{equation}

with the asymptotic behavior \(S_j(n) \sim (\ln n)^{-2} n^{-1}\) attained as \(\rho \to \rho_{max}\).

In the reference configuration defined by Equation \eqref{final_Sj(n)},
the parameters are fixed to $C_2=0.66$, $r=2/3$, and $\rho_{\max}=2$,
corresponding respectively to the Hardy-Littlewood twin prime constant,
the normalization implied by $S_{\ln n}\sim(\ln n)^{-2}$,
and the upper bound on normalized gaps suggested by the Cramér conjecture.

For moderate values of $n$ ($n<10^6$), the quantity $nS_j(n)$ significantly
underestimates $\mathrm{Li}(n)$.
This behavior is consistent with the fact that
the Hardy-Littlewood prediction $C_2(n)\to 0.66$ is asymptotic in nature,
while for finite ranges the empirical twin prime density
$K_2(n)$ can substantially exceed its limiting form
$2C_2(\ln n)^{-2}$ \cite{Granville1995}.
As $n$ increases to the range $10^{30}\lesssim n\lesssim 10^{80}$,
the agreement with the Prime Number Theorem constraint improves,
and the relative error $|\varepsilon(n)|$ remains below $2.5\%$.

However, for larger values of $n$ the discrepancy re-emerges:
for $n\gtrsim 10^{80}$ we observe $|\varepsilon(n)|\gtrsim 3\%$,
with a monotonic growth reaching $|\varepsilon(n)|>10\%$ for $n>10^{300}$.
This behavior indicates that the simultaneous enforcement of
$C_2=0.66$, $r=2/3$, and $\rho_{\max}=2$ is not asymptotically
compatible with the Prime Number Theorem within the present framework.

To identify the source of this misalignment between $nS_j(n)$ and $\mathrm{Li}(n)$,
we analyze three complementary scenarios, each relaxing one of the above assumptions.

\medskip
\noindent
\textit{(i) Variable twin-prime constant.}
Fixing $r=2/3$ and $\rho_{\max}=2$, we determine the value of $C_2(n)$
required to enforce the Prime Number Theorem.
For $n=10^{300}$ this yields $C_2(n)\simeq 1$ and $k_1(n)\simeq 68\%$.
As $n\to\infty$, one formally obtains $C_2(n)\to 1.5$ and $k_1(n)\to 1$,
which would imply an average prime gap smaller than $\ln n$,
in contradiction with the Prime Number Theorem.
This observation does not refute the Hardy-Littlewood conjecture,
but shows that within the present model $C_2(n)$ cannot grow indefinitely,
and that additional structural parameters must contribute to the shape of $S_j(n)$.

\medskip
\noindent
\textit{(ii) Variable normalization.}
Fixing $C_2=0.66$ and $\rho_{\max}=2$, we determine the normalization $r(n)$
required to satisfy the Prime Number Theorem.
Although this choice reduces the global error,
it forces $k_1(n)$ to increase with $n$,
violating the intrinsic relation \eqref{k1-C2}
between $k_1(n)$ and $C_2$.
This supports the stability of the condition
$S_{\ln n}\sim(\ln n)^{-2}$ and the associated value $r=2/3$.

\medskip
\noindent
\textit{(iii) Variable maximal gap scale.}
Fixing $C_2=0.66$ and $r=2/3$, we instead allow $\rho_{\max}$ to depend on $n$.
This scenario preserves both the Prime Number Theorem
and Equation \eqref{k1-C2},
but requires $\rho_{\max}(n)>2$ for sufficiently large $n$,
which is incompatible with the Cramér conjecture.
Within the present framework, this suggests that
the bound $j < (\ln n)^2$ cannot hold uniformly in $n$.

\medskip
\noindent
Table \ref{table:17} summarizes the behavior of $k_1(n)$ and $\varepsilon(n)$,
highlighting how the Prime Number Theorem constraint
can be restored by relaxing $C_2$, $r$, or $\rho_{\max}$ individually.

\begin{table}[h!]
\centering
\begin{tabular}{l c r c c c}
 \hline
 $ n $         & $k_1(n)$  & $\varepsilon(n)$   & $C_2(n)$  & $r(n)$ & $\rho_{max}(n)$    \\ [1ex]
 \hline
 $10^{10}$     & $55.6 \%$  & $2.11 \%$      & $0.620$   & $0.646$  & $1.955$  \\
 $10^{25}$     & $57.0 \%$  & $2.65 \%$      & $0.602$   & $0.642$  & $1.941$       \\
 $10^{50}$     & $58.9 \%$  & $0.11 \%$      & $0.658$   & $0.666$  & $1.998$    \\
 $10^{100}$    & $60.3 \%$  & $-3.66 \%$     & $0.761$   & $0.700$  & $2.095$      \\
 $10^{200}$    & $62.5 \%$  & $-7.95 \%$     & $0.903$   & $0.742$  & $2.230$    \\
 $10^{300}$    & $63.9 \%$  & $-10.51 \%$    & $1.002$   & $0.768$  & $2.325$  \\ [1ex]
 \hline
\end{tabular}
\caption{$k_1(n)$ and $\varepsilon(n)$ in the standard case, $C_2(n)$, $r(n)$ and $\rho_{max}(n)$ required to verify the Prime Number Theorem (as 3 different scenarios).}
\label{table:17}
\end{table}

Motivated by the admissible range $C_2(n)\in(0.5,1.5)$,
we propose the phenomenological form
\begin{equation}
C_2(n)\sim \frac{1}{2}+\frac{1}{1+\mu}
=\frac{1}{2}+\chi(n),
\label{C2}
\end{equation}
which provides a substantially improved fit to empirical data
over a wide range of $n$ (Table \ref{table:19}),
compared to the constant value $C_2=0.66$.

\begin{table}[h!]
\centering
\begin{tabular}{l c c r r}
 \hline
 $ n $       & $C_{emp}$ & $C_2(n)$  & $C_2(n)/C_{emp} - 1$  & $0.66/C_{emp} - 1$    \\ [1ex]
 \hline
 $10^{3}$   & $0.835$ & $0.841$  & $0.71 \%$ & $-20.94 \%$       \\
 $10^{6}$   & $0.780$ & $0.776$  & $-0.49 \%$  & $-15.32 \%$            \\
 $10^{9}$    & $0.735$ & $0.748$  & $1.73 \%$  & $-10.22 \%$         \\
 $10^{12}$    & $0.714$ & $0.732$  & $2.45 \%$  & $-7.55 \%$         \\
 $10^{15}$   & $0.702$ & $0.720$ & $2.56 \%$ & $-5.98 \%$            \\ [1ex]
 \hline
\end{tabular}
\caption{Empirical data of $C_2(n)$, expected $C_2(n)$ from Equation \eqref{C2}, errors with respect to empirical data for $C_2(n)$ and $C_2 \sim 0.66$.}
\label{table:19}
\end{table}

Within this parametrization, $C_2(n)$ becomes lower than $0.66$ for $n\gtrsim 10^{85}$,
indicating that the Hardy-Littlewood constant,
if interpreted as a strict asymptotic limit,
would not describe the observed scaling.
We emphasize that this conclusion is conditional on the present model
and does not constitute a proof or disproof of the conjecture itself.

Combining Equations \eqref{C2}, \eqref{Sjn1} and \eqref{Sjn2},
we obtain
$k_1(n)\sim \tfrac{1}{2}(1+\chi)$ and $k_2(n)\sim \tfrac{1}{2}(1-\chi)$,
with corresponding average densities
$\langle S_j(n)\rangle_1\sim(1+\chi)(\ln n)^{-2}$
and
$\langle S_j(n)\rangle_2\sim(1-\chi)(\ln n)^{-(\rho_{\max}+1)}$.
In the limit $n\to\infty$, this yields
$C_2(n)\to 0.5$ and $k_1(n)=k_2(n)=50\%$.

Finally, enforcing the Prime Number Theorem under the assumption \eqref{C2}
requires $\rho_{\max}(n)>2$ for $n\gtrsim 10^{61}$,
a regime beyond currently accessible empirical data.
A quadratic fit for $\rho_{\max}(n)$ is the following:

\begin{equation}
\rho_{max}(n) \sim 2.192 - 0.2257 \mu + 0.03828 \mu^2
\label{rho_max_fit}
\end{equation}

Table \ref{table:18} summarizes the behavior of the model under the most general configuration considered in this work and shows that relations \eqref{k1-C2} between $C_2(n)$ and $k_1(n)$ are also satisfied empirically.
In this setting, the twin-prime density is described by a scale-dependent coefficient $C_2(n)$ given by Equation \eqref{C2}, which provides an excellent match to available empirical data but does not converge to the Hardy-Littlewood constant.
The normalization is fixed to $r=2/3$, enforcing the scaling $S_{\ln n}\sim(\ln n)^{-2}$,
while $\rho_{\max}(n)$ is determined self-consistently by imposing the Prime Number Theorem as a global constraint, and therefore forcing $\rho_{\max}(n)$ to exceed $2$ for sufficiently large $n$.

Within the present framework, the simultaneous enforcement of
the Hardy-Littlewood constant, the Cramér upper bound,
and the Prime Number Theorem leads to an internal inconsistency.
This does not constitute evidence against any of these conjectures individually,
but highlights the rigidity imposed by combining them into a single
global density model.

\begin{table}[h!]
\centering
\begin{tabular}{l c c c r r}
 \hline
 $ n $     & $\mu$  & $C_2(n)$ & $\rho_{max}(n)$  & $k_1(n)$ & $\varepsilon(n)$   \\ [1ex]
 \hline
 $10^{10}$   & $3.137$ & $0.742$   & $1.861$  & $61.0 \%$   & $-0.10 \%$      \\
 $10^{25}$   & $4.053$ & $0.698$   & $1.906$  & $60.1 \%$  & $0.19 \%$           \\
 $10^{50}$   & $4.746$ & $0.674$   & $1.983$  & $59.6 \%$  & $-0.03 \%$        \\
 $10^{100}$  & $5.439$ & $0.655$   & $2.097$  & $58.0 \%$  & $-0.12 \%$          \\
 $10^{200}$  & $6.132$  & $0.640$  & $2.247$  & $56.9 \%$  & $-0.04 \%$    \\
 $10^{300}$  & $6.538$  & $0.633$  & $2.353$  & $56.3 \%$  & $0.04 \%$    \\ [1ex]
 \hline
\end{tabular}
\caption{$\mu$, $C_2(n)$ from Equation \eqref{C2}, $\rho_{max}(n)$ from Equation \eqref{rho_max_fit}, resulting $k_1(n)$ and $\varepsilon(n)$.}
\label{table:18}
\end{table}

\subsection{Quantitatve estimates on extreme gaps}

We now discuss the general implications of the proposed framework
for large and extreme prime gaps.

\medskip
\noindent
\textit{(i) Asymptotic decay and threshold behavior.}
For $\rho \sim 2$, Equation \eqref{final_Sj(n)_2} yields the upper-envelope behavior
\begin{equation}
S_j(n) \;\sim\; (\ln n)^{-2}\, n^{-(\rho-1)^3},
\end{equation}
which interpolates smoothly between the typical gap regime ($\rho \lesssim 1$)
and the extreme tail of the distribution.
Imposing the condition $S_j(n) \sim n^{-1}$,
corresponding to an expected multiplicity of order one,
leads to an approximate threshold value $\rho_T$ satisfying
\begin{equation}
\rho_T \sim 1 + \left( \frac{\ln n - 2\mu}{\ln n} \right)^{1/3},
\label{rho_T}
\end{equation}
with $\rho_T \to 2$ for $n \to \infty$.

This result identifies $\rho=2$ as a \emph{critical probabilistic threshold},
rather than as a rigid upper bound.
While prime gaps with $\rho \approx \rho_T$ are not guaranteed to occur for a given $n$,
since $nS_j(n)$ becomes subunitary,
the model does not exclude the existence of gaps with $\rho > \rho_T$.

In this sense, the scale $\rho=2$ emerges dynamically from the interaction
between the Prime Number Theorem constraint and the decay of $S_j(n)$,
without being imposed a priori.
This interpretation is consistent with classical heuristic arguments
underlying Cramér's conjecture \cite{Cramer1936},
while allowing for deviations in the extreme tail.

\medskip
\noindent
\textit{(ii) Super gaps and their expected multiplicity.}
Motivated by the above considerations, we define \emph{super gaps} as prime gaps satisfying $\rho > 2$.
We use the precise expression from Equation \eqref{final_Sj(n)_2} evaluated near the critical scale $\rho \sim 2$,
giving
\begin{equation}
S_{(\ln n)^2}(n)\sim \exp\!\left(-3\mu - \frac{\ln n - \mu}{(\rho_{\max}-1)^3}\right),
\end{equation}
For $\rho>2$ we have the trivial condition $S_j(n)<S_{(\ln n)^2}(n)$.

On this basis, an upper bound for the expected number of super gaps up to scale $n$ can be obtained via
\begin{equation}
N_S(n) < \frac{1}{2}\,n\,S_{(\ln n)^2}(n)\,\big[(\ln n)^{\rho_{\max}} - (\ln n)^2\big],
\end{equation}
In order to compare with the total number of prime gaps, we define the relative frequency
\begin{equation}
R_S(n) = N_S(n)\,\mathrm{Li}(n)^{-1},
\end{equation}
which measures the expected ratio between super gaps and all prime gaps.

Although $N_S(n)$ grows with $n$ due to the overall normalization of $S_j(n)$,
the ratio $R_S(n)$ decays rapidly with $n$, implying that super gaps remain extremely rare.
Numerical evaluation (Table \ref{table:9}) indicates that the first super gaps are only expected for extremely large scales, on the order of $n\sim10^{61}$,
at which the relative frequency $R_S(n)\sim10^{-57}$.

The present estimates are probabilistic in nature and address expected
frequencies rather than deterministic constructions.
They are therefore complementary to the explicit constructions of large gaps
developed in \cite{Westzynthius1931,Rankin1938,FGKMT2016},
which demonstrate existence but do not describe global statistical behavior.

From a complementary perspective,
only prime gaps with length $j < j_{\max}(n)$
can be expected to occur with non-negligible probability.
Writing
\begin{equation}
j_{\max}(n) \;=\; n^{\theta(n)},
\qquad
\theta(n) \;\sim\; \frac{\mu}{\ln n}\,\rho_{\max}(n),
\label{theta}
\end{equation}
we observe that $\theta(n) \to 0$ as $n \to \infty$,
even when $\rho_{\max}(n) > 2$ (Table \ref{table:9}).
Consequently, the maximal gap scale remains subpolynomial in $n$,
in agreement with all known rigorous results on large prime gaps
\cite{FGKMT2016}.

\begin{table}[h!]
\centering
\begin{tabular}{l c c c c c}
 \hline
 $ n $         & $\rho_{max}(n)$    & $j_{max}(n)$     & $\theta(n)$   & $N_S(n)$ & $R_S(n)$ \\ [1ex]
 \hline
 $10^{10}$    & $ 1.861 $  & $3.42\cdot 10^2$   & $0.253$       & $-$ & $-$ \\
 $10^{25}$     & $1.906$    & $2.26\cdot 10^3$     & $0.134$        & $-$ & $-$ \\
 $10^{50}$    & $1.983$     & $1.22\cdot 10^4$    & $0.082$        & $-$ & $-$ \\
 $10^{100}$    & $2.097$     & $8.98\cdot 10^4$   & $0.050$        & $2.42 \cdot 10^{21}$ & $5.54 \cdot 10^{-75}$ \\
 $10^{200}$     & $2.247$     & $9.67\cdot 10^5$   & $0.030$        & $2.13 \cdot 10^{93}$ & $9.77 \cdot 10^{-103}$  \\
 $10^{300}$   & $2.353$     & $4.78\cdot 10^6$   & $0.022$        & $4.77 \cdot 10^{174}$ & $3.29 \cdot 10^{-121}$  \\
 \hline
\end{tabular}
\caption{Expected upper limits for $\rho$ and $j$, the coefficient $\theta(n)$, the expected number of super gaps ($\rho>2$) and their expected relative frequency $R_S(n)$, as a function of $n$.}
\label{table:9}
\end{table}

\begin{figure}[ht]
\begin{center}
  \includegraphics[width=0.5\columnwidth]{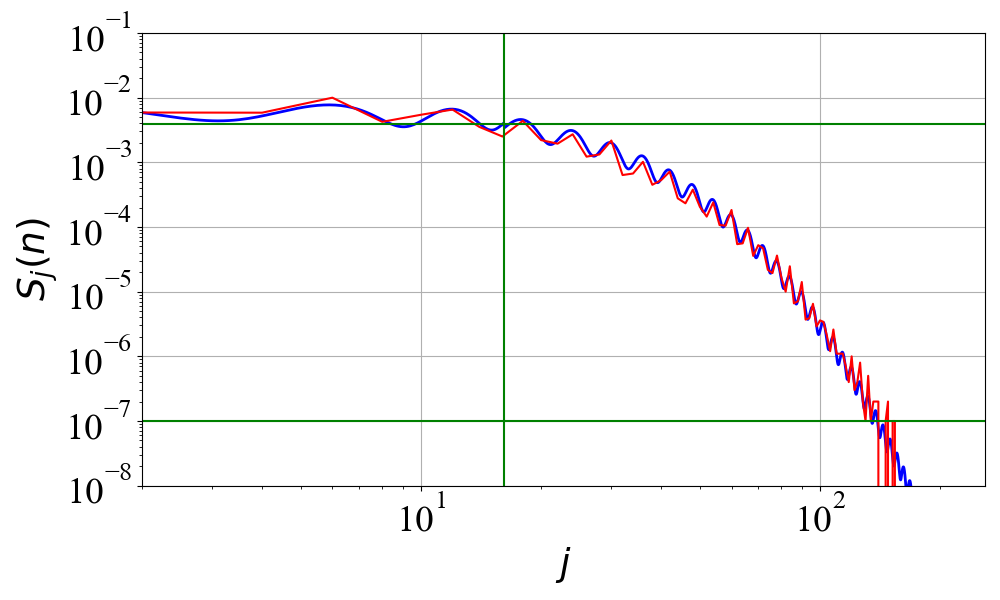}
  \caption{$S_j(n)$ \eqref{final_Sj(n)_2} (blue) compared to $P_{emp}(j,n)$ (red), for $n=10^7$ and $\rho_{max} = 1.92$. As references, the horizontal green lines are $S_j(n) = (\ln n)^{-2}$ and $S_j(n) = n^{-1}$, while the vertical green line is $j = \ln n$.}
  \label{S(j)_10_7}
\end{center}
\end{figure}

\begin{figure}[h!]
\begin{center}
  \includegraphics[width=0.5\columnwidth]{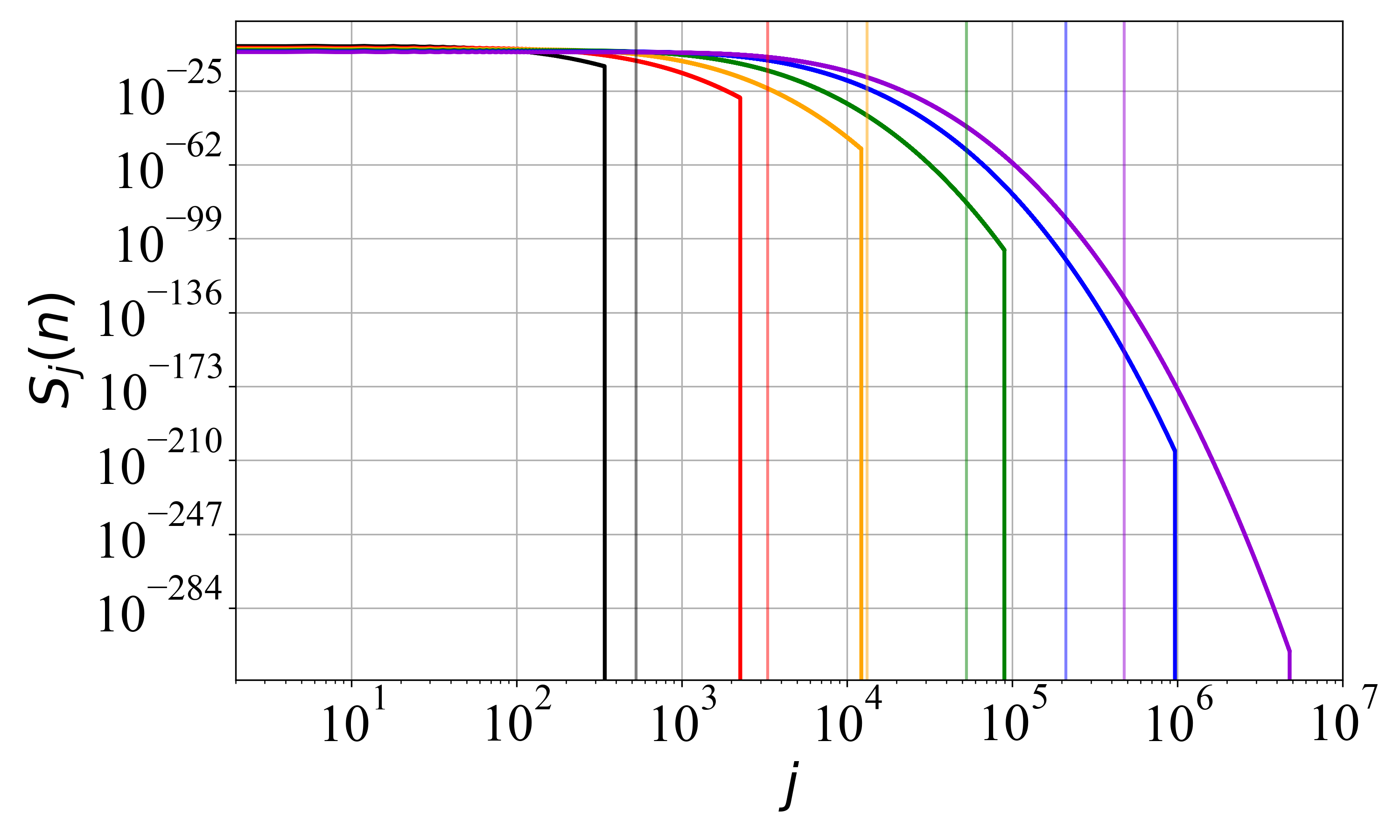}
  \caption{$S_j(n)$ \eqref{final_Sj(n)_2} as a function of $j$: $n = 10^{10}$ (black), $n = 10^{25}$ (red), $n = 10^{50}$ (orange), $n = 10^{100}$ (green), $n = 10^{200}$ (blue), $n = 10^{300}$ (violet). Vertical lines correspond to $\rho = 2$ for the considered values of $n$.}
  \label{S(j)_j}
\end{center}
\end{figure}

\begin{figure}[h!]
\begin{center}
  \includegraphics[width=0.5\columnwidth]{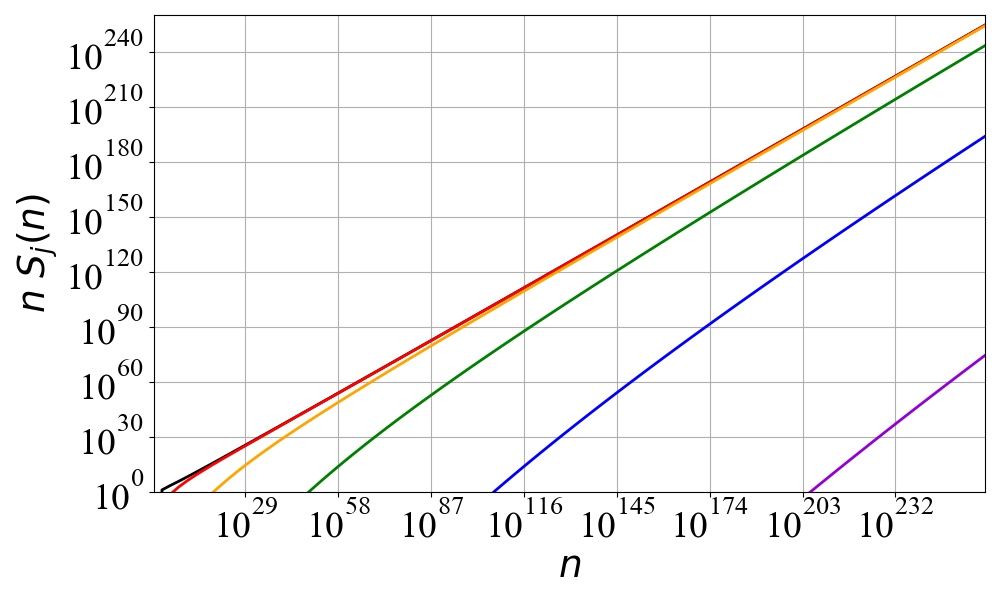}
  \caption{$nS_j(n)$ \eqref{final_Sj(n)_2} as a function of $n$: $j = 10$ (black), $j = 10^{2}$ (red), $j = 10^{3}$ (orange), $j = 10^{4}$ (green), $j = 10^{5}$ (blue), $j = 10^{6}$ (violet).}
  \label{S(j)_n}
\end{center}
\end{figure}

\medskip
\noindent
\textit{(iii) Comparison with empirical data and numerical extrapolations.}
Within this work, the density $S_j(n)$ is tested against empirical data
for all prime gaps $j \in (2,\lfloor (\ln n)^2 \rfloor)$
up to $n = 10^7$ (Figure \ref{S(j)_10_7}).
The quantity $nS_j(n)$ is evaluated numerically up to $n = 10^{300}$, far beyond the empirically accessible range, providing a controlled extrapolation of the model (Table \ref{table:8}).
Figures \ref{S(j)_j} and \ref{S(j)_n}
illustrate, respectively, the dependence of $S_j(n)$ on $j$
at fixed $n$ and the scaling of $nS_j(n)$ as a function of $n$
for representative gap sizes.

Overall, these results show that the proposed framework
reproduces the observed behavior of small and intermediate gaps,
while yielding nontrivial and testable predictions
for the extreme tail of the distribution.

In particular, the existence of large local deviations from average prime density, as demonstrated by Maier’s theorem on primes in short intervals \cite{Maier1985}, supports the need for a global density-based perspective rather than purely local probabilistic models.

Moreover, the emergence of a dynamic threshold near $\rho=2$,
the quantitative estimate of super-gap multiplicities,
and the subpolynomial growth of $j_{\max}(n)$
constitute robust and original features of the model,
largely independent of unproven conjectures.

\begin{table}[h!]
\centering
\begin{tabular}{l l l l l l}
 \hline
 $ \ $             & $10^2$      & $10^3$         & $10^4$    & $10^5$   \\ [1ex]
 \hline
 $10^{25}$        & $1.11\cdot 10^{21}$    & $8.43\cdot 10^{8}$    & $-$  & $-$ \\
 $10^{50}$         & $6.97\cdot 10^{45}$    & $7.09\cdot 10^{39}$   & $5.67\cdot 10^{1}$  & $-$ \\
 $10^{100}$         & $1.83\cdot 10^{95}$    & $9.95\cdot 10^{92}$    & $4.09\cdot 10^{68}$ & $-$  \\ 
 $10^{200}$        & $4.69\cdot 10^{194}$   & $1.07\cdot 10^{194}$    & $1.55\cdot 10^{180}$  & $2.72\cdot 10^{123}$ \\ 
 $10^{300}$         & $2.11\cdot 10^{294}$   & $1.14\cdot 10^{294}$   & $5.32\cdot 10^{284}$  & $8.95\cdot 10^{238}$  \\ [1ex]
 \hline
\end{tabular}
\caption{$nS_j(n)$ from Equation \eqref{final_Sj(n)_2}, as a function of $j$ (columns) and $n$ (rows), considering $C_2(n)$ \eqref{C2} and $\rho_{max}$ \eqref{rho_max_fit}.}
\label{table:8}
\end{table}

\newpage

\section{Applications to Goldbach-type problems}
\label{section 3}

The Goldbach conjecture is one of the most prominent open problems in analytic number theory.
It asserts that every even integer $2n \geq 4$ can be written as the sum of two prime numbers.
Equivalently, writing such a representation in symmetric form,
one may express it as
\begin{equation}
n = \frac{p_1 + p_2}{2} = p_1 + j = p_2 - j ,
\end{equation}
for some integer $j \geq 0$.
Thus, the conjecture can be reformulated by stating that every integer $n \geq 2$
is either itself prime or lies at the center of at least one pair of primes
separated by a distance $2j$.

Following the probabilistic and density-based viewpoint
introduced by Hardy and Littlewood in their classical work on additive problems
\cite{HW1923},
we recast the Goldbach conjecture in terms of the prime indicator function
$I(n)$ defined in Section \ref{section 1}.
In this language, the conjecture can be written as
\begin{equation}
\exists \, j \in [0,n-2] \;:\;
I(n-j)\, I(n+j) = 1
\qquad \forall\, n \geq 2 .
\end{equation}

Given the constraints $p_1 \geq 2$ and $p_2 \leq 2n-2$,
we define the number of prime pairs centered at $n$
(allowing the degenerate case $p_1 = p_2 = n$)
as
\begin{equation}
G(n)
\;=\;
\sum_{j=0}^{n-2} I(n-j)\, I(n+j) .
\label{G(n)}
\end{equation}
The Goldbach conjecture is therefore equivalent to the assertion that
$G(n) \geq 1$ for all $n \geq 2$.

Using the identity $I(n) = K(n) - K(n-1)$,
where $K(n)$ denotes the cumulative prime counting function introduced earlier,
Equation \eqref{G(n)} can be rewritten as
\begin{equation}
G(n)
=
\sum_{j=0}^{n-2}
\bigl[K(n-j)-K(n-j-1)\bigr]
\bigl[K(n+j)-K(n+j-1)\bigr] .
\end{equation}

In this formulation, the summand can be interpreted,
at the level of expected values, as encoding the local probability 
that both integers $n-j$ and $n+j$ are prime within a density-based model.
This reformulation shifts the focus of the problem
from the precise location of individual primes
to their local statistical distribution,
in the spirit of classical heuristic approaches
\cite{MontgomerySoundararajan2004}.

To illustrate this viewpoint,
consider an extreme hypothetical scenario in which all odd integers are prime.
In this case, the probability that both $n-j$ and $n+j$ are prime would be $1/2$,
and one would expect $G(n) \sim n/2$ representations of $2n$
as a sum of two primes.
More generally, even in the presence of local fluctuations or regions of low prime density,
it is the overall distribution of primes in intervals around $n$
that governs the typical size of $G(n)$.

Using the local density approximation introduced in Equation \eqref{D-L},
we write
\begin{equation}
K(n \pm j)
\;\approx\;
K(n) \pm D_j(n)\, j ,
\end{equation}
where $D_j(n)$ denotes the average prime density at distance $j$ from $n$.
At the level of expected values and to leading order,
this yields the approximation
\begin{equation}
G(n)
\;\sim\;
\sum_{j=0}^{n-2} D_j(n)^2 .
\label{G(n)_D(n)}
\end{equation}

Equation \eqref{G(n)_D(n)} should be understood as a heuristic relation
describing the expected number of Goldbach representations of $2n$
in terms of the squared local prime density.
This perspective is consistent with the classical Hardy-Littlewood heuristics
and complements both rigorous partial results
and extensive computational verifications
\cite{OliveiraSilva2014}.

\subsection{Bounds on local density fluctuations}

The novelty of the present approach does not lie in introducing a new
Goldbach counting function, but in isolating structural density constraints
that force the expected size of $G(n)$ to grow well beyond the minimal
threshold required by the conjecture.
In this sense, Goldbach-type representations emerge as a low-complexity
additive phenomenon compared to the global constraints already imposed
by prime density.

We analyze the Goldbach problem through the quantity $G(n)$
by separating the contributions coming from the intervals $(2,n)$ and $(n,2n)$,
with the aim of deriving theoretical upper and lower bounds
for its expected magnitude.

\medskip
\noindent
\textit{Upper bound.}
In the case of a maximally efficient pairing scenario,
we assume that for every prime $p_2 \in (n,2n)$
there exists a corresponding prime $p_1 \in (2,n)$
such that $p_1 + p_2 = 2n$.
This configuration should not be interpreted as attainable,
but rather as an extremal benchmark compatible with global density constraints.

Under this idealized configuration,
each prime in $(n,2n)$ contributes at most one Goldbach representation,
yielding the theoretical upper bound
\begin{equation}
G_{\max}(n)
\;\sim\;
\text{Li}(2n) - \text{Li}(n)
\;\sim\;
n \,
\frac{2\ln n - \ln(2n)}{\ln n \, \ln(2n)} .
\label{Gmax(n)}
\end{equation}
This bound corresponds to the maximal number of admissible centered pairs
compatible with the Prime Number Theorem
and serves as a natural ceiling for $G(n)$.

\medskip
\noindent
\textit{Lower bound.}
In the opposite extreme, we consider a minimally symmetric scenario,
in which the existence of a prime $p_1 < n$
complementing a given prime $p_2 > n$
is governed solely by the local prime density.
Let $j_2$ denote the average prime gap in the interval $(n,2n)$.
From the Prime Number Theorem, we have
\begin{equation}
j_2
\;\sim\;
\frac{\ln n \, \ln(2n)}{2\ln n - \ln(2n)}
\;>\;
\ln(2n) .
\end{equation}
In this setting, the probability that a given $p_2$
admits a complementary prime $p_1$
is of order $j_2^{-1}$.

Allowing this probability to be significantly smaller would,
within a density-based interpretation,
force a macroscopic imbalance between the prime counts in $(2,n)$ and $(n,2n)$,
in contradiction with the Prime Number Theorem.

Substituting the estimate $D_j(n) \sim j_2^{-1}$
into Equation \eqref{G(n)_D(n)},
we obtain the following lower bound
for the expected value of $G(n)$:
\begin{equation}
G_{\min}(n)
\;\sim\;
\frac{\bigl(\text{Li}(2n)-\text{Li}(n)\bigr)^2}{n}
\;\sim\;
n \left(
\frac{2\ln n - \ln(2n)}{\ln n \, \ln(2n)}
\right)^2 .
\label{Gmin(n)}
\end{equation}

\medskip
\noindent
\textit{Implications.}
Combining the above estimates,
we obtain the heuristic inequality
\begin{equation}
G_{\min}(n)
\;\lesssim\;
\mathbb{E}[G(n)]
\;\lesssim\;
G_{\max}(n) ,
\label{G_constraint}
\end{equation}
which constrains the expected number of Goldbach representations
within a narrow range determined solely by prime densities.
In particular, the condition $G_{\min}(n) \gg 1$
is significantly stronger than the requirement $G(n) \geq 1$,
and suggests that violations of the Goldbach conjecture
would require an anomalously sparse prime distribution.

To appreciate the rigidity of this bound,
it is instructive to examine the hypothetical threshold
at which the expected number of representations becomes marginal.
Setting $G_{\min}(n) \sim 1$
would formally correspond to an average prime density
$D(n) \sim n^{-1/2}$.
Such a decay is incompatible both with the Prime Number Theorem,
which implies $n/\ln n \gg \sqrt{n}$,
and with the gap statistics described by the $S_j(n)$ formalism in Section \ref{section 2},
which predicts that gaps larger than $\sqrt{n}$
cannot occur beyond moderate values of $n$.
Consequently, within this probabilistic framework,
the mechanisms required to falsify the Goldbach conjecture
appear incompatible with the known global behavior of primes.

\begin{figure}[ht]
\begin{center}
  \includegraphics[width=0.5\columnwidth]{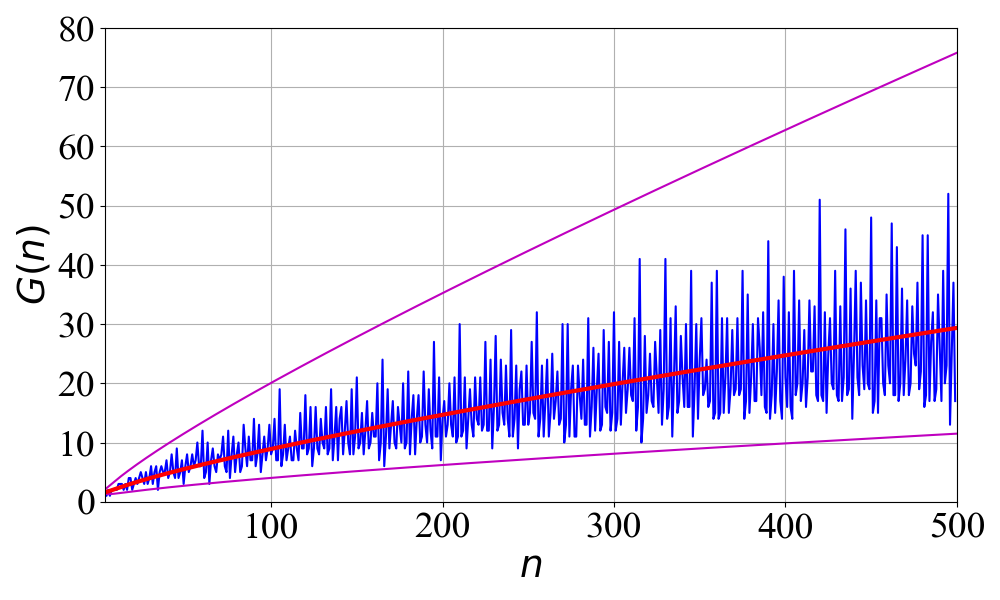}
  \caption{$G(n)$ \eqref{G(n)} (blue), $\langle G(n) \rangle$ \eqref{G average} (red), $G_{max}(n)$ \eqref{Gmax(n)} and $G_{min}(n)$ \eqref{Gmin(n)} (magenta), for $n \leq 500$.}
  \label{G(n)_fig}
\end{center}
\end{figure}

Finally, averaging over admissible configurations,
the typical size of $G(n)$ is well approximated by
\begin{equation}
\langle G(n) \rangle
\;\sim\;
\frac{n}{\ln(3n/2)^{3/2}} ,
\label{G average}
\end{equation}
in good agreement with numerical evaluations
for moderate values of $n$
(Figure~\ref{G(n)_fig}).

\subsection{Compatibility with Hardy–Littlewood heuristics}

The classical Hardy-Littlewood conjecture
provides a refined asymptotic prediction for the Goldbach counting function,
asserting that
\begin{equation}
G(n)
\;\sim\;
2 C_2 \,
\frac{n}{(\ln n)^2}
\prod_{\substack{p \mid n \\ p>2}}
\frac{p-1}{p-2},
\label{HLB}
\end{equation}
where $C_2$ denotes the twin prime constant
and the finite Euler product accounts for local congruence obstructions
arising from residue classes modulo small primes.

The present framework does not assume this conjecture,
nor does it introduce arithmetic correction factors \emph{a priori}.
Instead, it derives upper and lower bounds for $G(n)$
directly from global prime densities and gap statistics.
The bounds obtained in
Equations ~\eqref{Gmax(n)} and \eqref{Gmin(n)}
are therefore model-internal and non-conjectural within the adopted statistical formalism.
Remarkably, they are fully compatible with the Hardy-Littlewood prediction
at the level of order of magnitude.

This compatibility can be made more explicit.
Empirically, for $n$ in the range $10^1 \le n \le 10^4$,
the product over primes dividing $n$ satisfies
\begin{equation}
1
\;<\;
\prod_{\substack{p \mid n \\ p>2}}
\frac{p-1}{p-2}
\;\lesssim\;
1.4 \,\ln n,
\label{HL_product_bound}
\end{equation}
with the effective coefficient of $\ln n$ decreasing as $n$ increases.
In particular, this product grows at most logarithmically
and remains uniformly bounded away from zero.

Moreover, within the present framework,
the twin-prime constant $C_2$ may be consistently replaced
by the scale-dependent quantity $C_2(n)$.
Inserting Equation\eqref{C2} into the Hardy-Littlewood expression \eqref{HLB},
and using the empirical bound \eqref{HL_product_bound},
we obtain the inequalities
\begin{equation}
\frac{n}{(\ln n)^2}
\lesssim
G(n) \left(\frac{1 + \mu}{3+\mu} \right)
\lesssim\
\frac{1.4\,n}{\ln n}.
\label{HL_bounds_mu}
\end{equation}

In the limit $n \to \infty$,
both bounds asymptotically reduce to $n/(\ln n)^2 \lesssim G(n) \lesssim n/\ln n$, 
which exactly reproduces the upper and lower limits
previously derived in a completely independent way
within the prime-gap formalism.

This observation has an important conceptual consequence.
The lower bound obtained here is far stronger than the mere condition
$G(n) \geq 1$ required by the strong Goldbach conjecture.
Rather than relying on specific arithmetic correlations,
the existence of Goldbach representations emerges
as a structural consequence of prime density alone.

\subsection{Structural and probabilistic insights}

From a combinatorial perspective,
it is natural to consider the cumulative number of prime pairs
centered at integers up to a given scale.
We define
\begin{equation}
G_*(n) \;=\; \sum_{j=2}^{n} G(j),
\end{equation}
which counts the total number of representations
of even integers $2j \leq 2n$
as sums of two primes, allowing multiplicities.

Independently of Goldbach-type constraints,
the total number of unordered prime pairs
$(p_1,p_2)$ with $p_1,p_2 \in [2,n]$
and $p_1 \leq p_2$
is given by
\begin{equation}
C_P(n)
\;=\;
K(n)\,
\left\lfloor \frac{K(n)+1}{2} \right\rfloor,
\label{C_P_1}
\end{equation}
where $C_P(n) \to K(n)^2/2$ for $n \to \infty$.

Indeed, $C_P(n)$ does not account for prime pairs
with $p_1 \leq n < p_2$,
while $C_P(2n-2)$ includes also pairs with $p_1 > n$
that cannot contribute to representations of integers $\leq 2n$.
Thus,
\begin{equation}
C_P(n)
\;\leq\;
G_*(n)
\;\leq\;
C_P(2n-2).
\end{equation}

However, distinct prime pairs may correspond to the same central integer.
That is, different pairs $(p_1,p_2)$ may satisfy
$p_1+p_2=2j$ for the same $j$.

We denote by $G_D(n)$ the number of such \emph{duplicate pairs}.
Accordingly, we define
\begin{equation}
G_G(n)
\;=\;
G_*(n) - G_D(n),
\end{equation}
which counts the number of integers $j \leq n$
that are either prime themselves
or lie at the center of at least one prime pair.

The purpose of this cumulative viewpoint is not to refine asymptotics,
but to assess whether purely combinatorial obstructions could plausibly
prevent full coverage of admissible centers.

The strong Goldbach conjecture is equivalent to the statement $G_G(n) = n-1 \ \forall \ n \geq 2$, 
namely that every integer $2 \leq j \leq n$
admits at least one Goldbach representation.

Trivially, one always has $G_G(n) < n$ and $G_D(n) \geq C_P(n) - n$,
since at least $C_P(n)-n$ pairs must map onto already covered centers.
This indicates that, even in the absence of duplicate pairs ($G_D(n) = 0$), we have $G_G(n) < K(2n)^2/2$ for $n \to \infty$.

A necessary, though not sufficient, condition
for the validity of the Goldbach conjecture
is therefore
\begin{equation}
G_G(n) = n -1 \Rightarrow K(n) > \sqrt{2n},
\end{equation}
which ensures that the number of distinct prime pairs
is large enough to potentially cover all centers.
This condition is asymptotically satisfied,
since $K(n) \sim n/\ln n \gg \sqrt{n}$.

Although purely combinatorial and far from sufficient,
this argument reinforces a recurring theme of this section:
any counterexample to Goldbach would require a breakdown
of prime density at a scale incompatible with both analytic
and probabilistic constraints.

\newpage

\section{Oscillatory analysis of Zeta stability and zero localization}
\label{section 4}

The Riemann zeta function $\zeta(s)$ is classically defined,
for complex $s$ with $\Re(s)>1$, by the absolutely convergent Dirichlet series
\begin{equation}
\zeta(s) \;=\; \sum_{n=1}^{\infty} \frac{1}{n^{s}} .
\end{equation}
It admits an analytic continuation to the complex plane,
with a unique simple pole at $s=1$,
and can be equivalently expressed through Euler's product formula
\begin{equation}
\zeta(s)
\;=\;
\prod_{p}
\frac{1}{1-p^{-s}},
\label{euler_product}
\end{equation}
where the product extends over all prime numbers.
This identity encodes the fundamental relation
between $\zeta(s)$ and the distribution of primes.

Within the density-based viewpoint developed in Section \ref{section 1},
one may interpret the divisibility condition $i \mid n$
in terms of its average frequency over integers,
which is asymptotically $i^{-1}$.
Accordingly, the complementary event occurs with average density $1-i^{-1}$.

Restricting attention to prime divisors only,
as in the Sieve of Eratosthenes,
avoids systematic overcounting
and leads to a consistent approximation of prime density
\cite{Tenenbaum1995}.

We consider a truncated version of the classical Mertens product,
which asymptotic behaviour is governed by Mertens' third theorem.
\begin{equation}
\Pi(n)
\;=\;
\prod_{p \leq \lfloor n/2 \rfloor}
\bigl(1-p^{-1}\bigr),
\end{equation}
which yields
\begin{equation}
K(n)
\;\sim\;
\sum_{j=2}^{n} \Pi(j)
\;\sim\;
\mathrm{Li}(n),
\end{equation}
in agreement with the Prime Number Theorem.

We define the truncated zeta function
\begin{equation}
\zeta(s,n)
\;=\;
\prod_{p \leq \lfloor n/2 \rfloor}
\frac{1}{1-p^{-s}},
\end{equation}

At a heuristic level, for $s = 1$ we have
\begin{equation}
\zeta(1,n) \sim \frac{n}{\mathrm{Li}(n)} \sim \ln n,
\end{equation}
which coincides with the expected average prime gap.

The divergence of $\zeta(s)$ at $s=1$
thus reflects the unbounded growth of prime gaps
in the limit $n \to \infty$,
a phenomenon consistent with both the Prime Number Theorem
and known results on large gaps between primes.

For $\Re(s)>1$, we define
\begin{equation}
M(s)
\;=\;
1-\zeta(s)^{-1},
\label{M(s)}
\end{equation}
which may be interpreted, in an average-density sense,
as the proportion of integers divisible by at least one prime power $p^s$.

The pole of $\zeta(s)$ at $s=1$ reflects the fact that the cumulative
exclusion of integers divisible by primes fails to converge,
corresponding to the absence of integers
free of prime divisors, in accordance with the Fundamental Theorem of Arithmetic.

As $\Re(s)$ increases,
the contribution of higher prime powers becomes increasingly sparse.
For instance, $\zeta(2)=\pi^{2}/6$
implies that the density of integers divisible by some $p^{2}$
is $1-6/\pi^{2}$,
with analogous behavior for larger values of $s$
(Table~\ref{table:21}).

\begin{table}[h!]
\centering
\begin{tabular}{c c c c c c c}
 \hline
 $ s $                        & $2$     & $3$   & $4$  & $5$  & $6$ & $7$  \\ [1ex]
 \hline
 $\zeta(s)$                 & $\pi^2/6$    & $1.202$ & $\pi^4/90$ & $1.037$  & $\pi^6/945$ & $1.008$ \\
 $M(s)$        & $39.21 \%$    & $16.81 \%$ & $7.61 \%$ & $3.57 \%$ & $1.70 \%$ & $0.79 \%$ \\ [1ex]
 \hline
\end{tabular}
\caption{$\zeta(s)$ and the resulting $M(s)$ \eqref{M(s)}, for $\Re(s) \in (2,7)$.}
\label{table:21}
\end{table}

For $\Re(s)\gg1$ and $\Im(s) = 0$,
one may write the first-order approximation
\begin{equation}
M(s) \sim \sum_{p=2}^{\lfloor n^{1/|s|} \rfloor} p^{-s}
\end{equation}
so that
\begin{equation}
\zeta(s) \sim \lim_{n \to \infty} \bigg( 1 - \sum_{p=2}^{\lfloor n^{1/|s|} \rfloor} p^{-s} \bigg)^{-1},
\end{equation}
reflecting the dominance of the leading sieve term
in the regime of rapidly decaying prime contributions.

\subsection{Representations of the Riemann Zeta function}

For $\Re(s)>1$, the Riemann zeta function admits several classical
representations in terms of fundamental arithmetic functions.
In particular, it is related to the von Mangoldt function
$\Lambda(n)$, defined by
\begin{equation}
\Lambda(n)=
\begin{cases}
\ln p, & \text{if } n=p^k \text{ for some prime } p, \\
0, & \text{otherwise}.
\end{cases}
\end{equation}
For $\Re(s)>1$, we have
\begin{equation}
\ln \zeta(s)
=
\sum_{n=2}^{\infty}
\frac{\Lambda(n)}{n^s \ln n},
\end{equation}
which makes explicit the contribution of prime powers
to the logarithmic structure of $\zeta(s)$.

Using Equation \eqref{I}, one may formally rewrite this expansion as
\begin{equation}
\ln \zeta(s)
\sim
\sum_{n=1}^{\infty}
\sum_{k=1}^{\infty}
\frac{\mathrm{Li}(n^{1/k})-\mathrm{Li}(n^{1/k}-1)}{k\,n^s},
\end{equation}
thereby expressing $\ln \zeta(s)$ in terms of local prime densities.

For $\Re(s)>0$, $\zeta(s)$ is related to the Dirichlet eta function
\begin{equation}
\eta(s) = \sum_{n=1}^{\infty} \frac{(-1)^{n-1}}{n^s} = (1-2^{1-s})\,\zeta(s),
\label{zeta-eta}
\end{equation}
which provides an analytic continuation of $\zeta(s)$ to the half-plane
$\Re(s)>0$, excluding the pole at $s=1$.

The zeta function satisfies the classical functional equation
\begin{equation}
\zeta(s)
=
2^s \pi^{s-1}
\sin\!\left(\frac{\pi s}{2}\right)
\Gamma(1-s)\,
\zeta(1-s),
\label{zeta-gamma}
\end{equation}
where $\Gamma(s)$ denotes the Euler gamma function.
Its Weierstrass product representation is given by
\begin{equation}
\Gamma(s)
=
\frac{e^{-\gamma s}}{s}
\prod_{n=1}^{\infty}
\left(1+\frac{s}{n}\right)^{-1}
e^{s/n},
\label{gamma(s)}
\end{equation}
with $\gamma \sim 0.57722$ denoting the Euler-Mascheroni constant,
\begin{equation}
\gamma = \lim_{n \to \infty} \bigg( \sum_{k=1}^{n} \frac{1}{k} - \ln n \bigg) = \ln \bigg( \prod_{n=1}^{\infty} \frac{e^{1/n}}{1 + 1/n}  \bigg).
\end{equation}

Since $\Gamma(1/2) = \sqrt{\pi}$, using Equation \eqref{gamma(s)}, $\gamma$ can be expressed as
\begin{equation}
\gamma = \ln \bigg( \frac{4}{\pi} \prod_{n=1}^{\infty} \frac{4n^2 e^{1/n}}{(2n+1)^2} \bigg)
\label{em_gamma}
\end{equation}
determining an expression of $\pi$ similar to the Wallis product:
\begin{equation}
\pi = 4 \prod_{n=1}^{\infty} \frac{n+1}{n+1+1/4n}
\label{greek pi}
\end{equation}

The Riemann zeta function is classically related to the Bernoulli numbers
$\{B_s\}$, defined recursively by
\begin{equation}
B_s
=
-\frac{1}{s+1}
\sum_{j=0}^{s-1}
\binom{s+1}{j}
B_j,
\end{equation}
and satisfying the well-known identity
\begin{equation}
\zeta(2s)
=
\frac{|B_{2s}|\, (2\pi)^{2s}}{2 (2s)!},
\label{zeta-bernoulli}
\end{equation}
for integers $s \geq 1$ \cite{Apostol}.

Within the present formalism, this relation suggests a natural
reinterpretation of Bernoulli numbers in terms of prime-based products.
Formally, one may write
\begin{equation}
B_s
=
-\frac{s!}{2^{\,s-1}\,\pi^s}
\cos\!\left(\frac{\pi s}{2}\right)
\prod_{n=1}^{\infty}
\frac{1}{1-I(n)\,n^{-s}},
\label{B_s}
\end{equation}
where $I(n)$ denotes the prime indicator function introduced in
Section \ref{section 1}.

This expression should be understood purely as a formal reconstruction,
obtained by combining Euler products with the classical identity
\eqref{zeta-bernoulli}, and is not intended as a new independent
definition of the Bernoulli numbers.

\begin{figure}[ht]
\begin{center}
  \includegraphics[width=0.5\columnwidth]{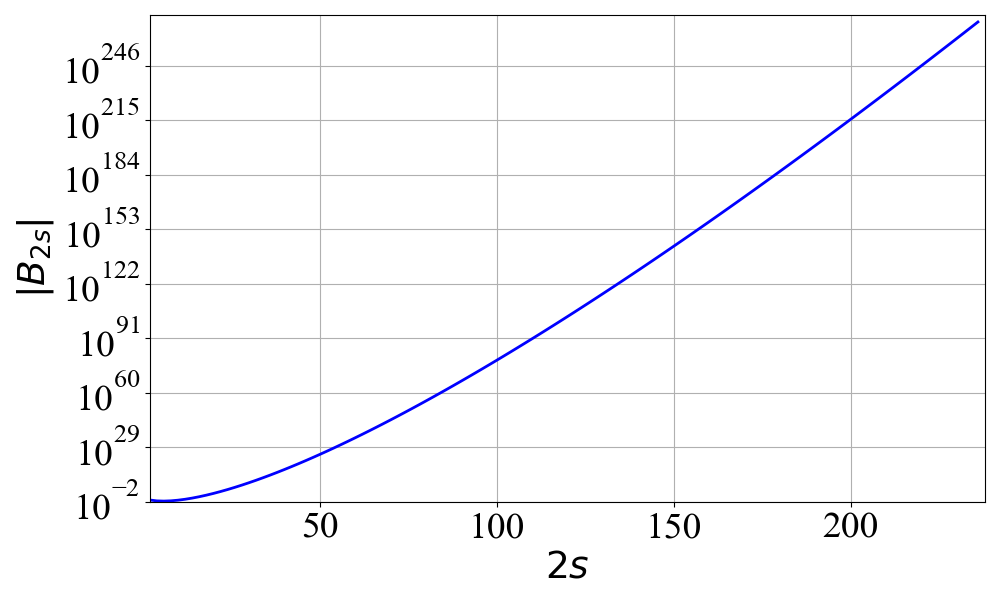}
  \caption{$|B_{2s}|$ \eqref{B_s_approx}, for $\Re(s) \leq 236$.}
  \label{B_2s}
\end{center}
\end{figure}

For $\Re(s)\gg1$, since $\zeta(s)\to1$ and using Stirling's approximation
$s!\sim\sqrt{2\pi s}\,(s/e)^s$, Equation \eqref{B_s} yields the classical
asymptotic behaviour
\begin{equation}
B_s
\sim
-\sqrt{8\pi s}
\left(\frac{s}{2\pi e}\right)^s
\cos\!\left(\frac{\pi s}{2}\right),
\label{B_s_approx}
\end{equation}
in agreement with standard results for Bernoulli numbers (Figure \ref{B_2s}).

The analytic structure of $\zeta(s)$ is further encoded in its
Hadamard product representation,
\begin{equation}
\zeta(s)
=
\frac{(2\pi)^s e^{s(\gamma/2-1)}}
{2(s-1)\Gamma(1+s/2)}
\prod_{\rho}
\left(1-\frac{s}{\rho}\right)
e^{s/\rho},
\label{hadamard}
\end{equation}
where the product runs over the non-trivial zeros $\rho$ of $\zeta(s)$
\cite{Titchmarsh}.

By substituting the explicit product representations for
$\Gamma(s)$, $\gamma$, and $\pi$ derived earlier,
one may formally rewrite $\zeta(s)$ in the factorized form
\begin{equation}
\zeta(s)
=
f(s)
\prod_{n=1}^{\infty}
g(s,n)\,h(n)^s\,e^{s/2n}
\prod_{\rho}
\frac{1-s/\rho}{e^{-s/\rho}},
\label{zeta(s)_new}
\end{equation}
where
\begin{equation}
f(s) = \frac{1}{4}\,\frac{s+2}{s-1}\left(\frac{4}{e}\right)^s, \qquad
g(s,n) = 1+\frac{s}{2(n+1)}, \qquad
h(n) = \frac{4n^3}{(4n^2-1)(n+1)}.
\end{equation}

In this representation, the pole at $s=1$ arises from the divergence of
$f(s)$, the trivial zeros correspond to the vanishing of $f(s)$ and
$g(s,n)$ at negative even integers, the normalization $h(n) \to 1$ for $n \to \infty$, and the non-trivial zeros are encoded
in the final product over $\rho$.
While algebraically equivalent to the classical Hadamard factorization,
this formulation highlights, at a heuristic level, the main analytic
features of $\zeta(s)$.

\subsection{Functional-equation heuristics for the critical line}

The celebrated Riemann Hypothesis (RH) asserts that all non-trivial zeros
of the Riemann zeta function lie on the critical line
$\Re(s)=1/2$.
This conjecture is known to be equivalent to sharp bounds on the error
term in the Prime Number Theorem and underlies much of modern analytic
number theory.

In the context of the symmetry-based framework developed in this work,
the functional equation of $\zeta(s)$ provides a natural starting point
for heuristic considerations about the location of its non-trivial zeros.

Let $\rho_1=a_1+ib_1$ and $\rho_2=1-a_1+ib_1$ denote a pair of non-trivial
zeros symmetric with respect to the critical line, with $0<a_1<1$.
The functional equation of $\zeta(s)$,
Equation \eqref{zeta-gamma},
naturally suggests considering the ratio
\begin{equation}
V(s)=\frac{\zeta(s)}{\zeta(1-s)},
\end{equation}
which explicitly encodes the functional symmetry
$s\leftrightarrow1-s$ at the level of analytic ratios.

Using the functional equation of the Gamma function,
\begin{equation}
\Gamma(s)=\frac{\pi}{\sin(\pi s)\,\Gamma(1-s)},
\end{equation}
together with $\sin(\pi s)=2\sin(\pi s/2)\cos(\pi s/2)$,
one obtains the explicit expression
\begin{equation}
V(s)=\frac{(2\pi)^s\,\Gamma(s)}{2\cos(\pi s/2)}.
\end{equation}

This expression shows that $V(s)$ is entirely determined by elementary
transcendental factors, independently of the zero structure of
$\zeta(s)$.
For purely exploratory and visualization purposes, we introduce the
auxiliary quantities
\begin{equation}
V_1(s)=\big|\ |V(s)|-1\ \big|^{k},
\qquad
V_2(s)=\left|\frac{\arg V(s)+\pi}{2\pi}\right|^{k},
\label{V1-V2}
\end{equation}
where $k>0$ is a visualization parameter with no theoretical significance.

The condition $|V(s)|\approx1$ reflects an approximate balance between the values of $\zeta(s)$ and
$\zeta(1-s)$ induced by the functional equation, while $\arg V(s)\approx-\pi$ corresponds to the phase behaviour
near non-trivial zeros.

\begin{figure}[ht]
\begin{center}
  \includegraphics[width=0.5\columnwidth]{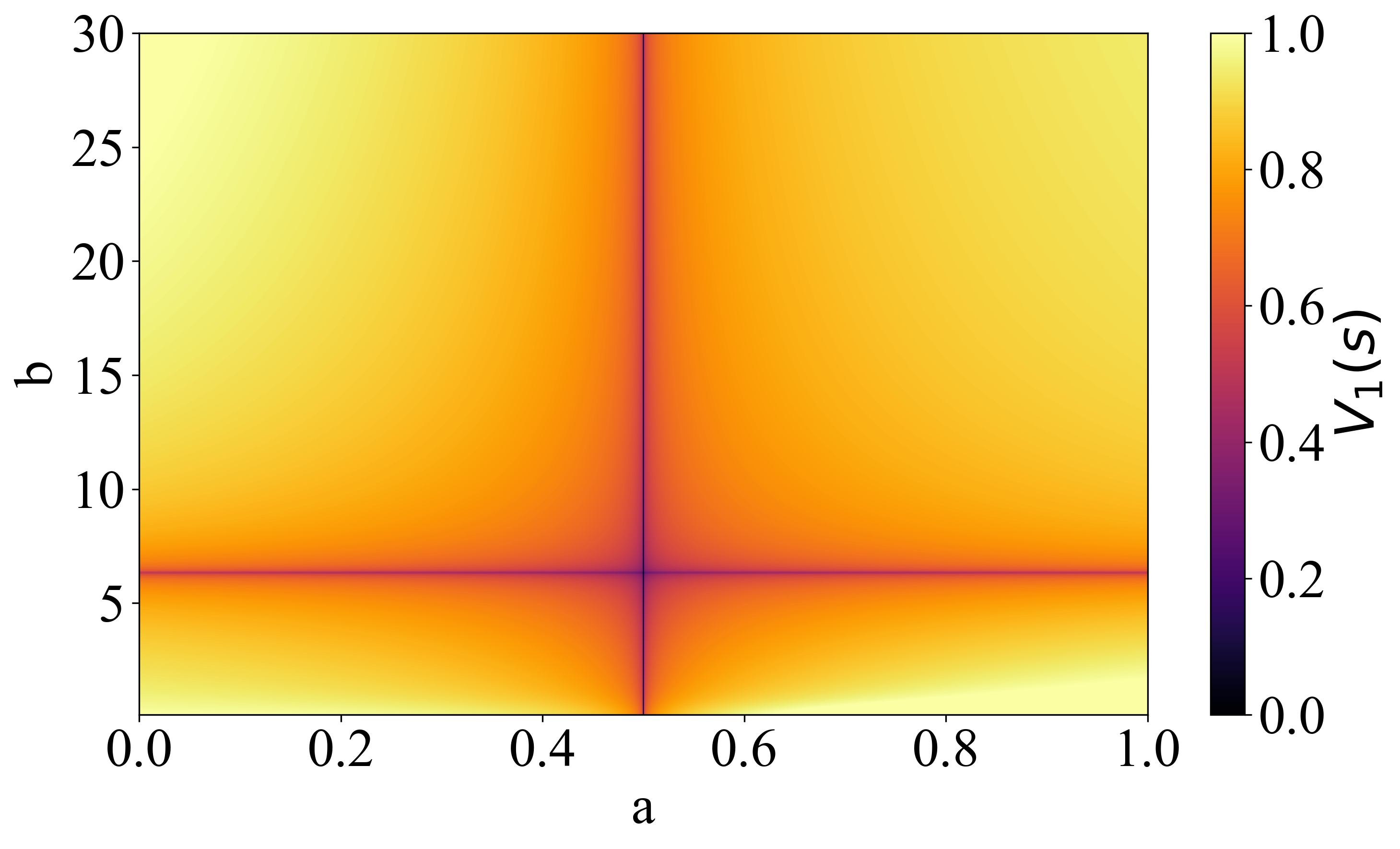}
  \caption{$V_1(s)$ \eqref{V1-V2} as a function of $a$ and $b$.}
  \label{V_1(s)}
\end{center}
\end{figure}

\begin{figure}[ht]
\begin{center}
  \includegraphics[width=0.5\columnwidth]{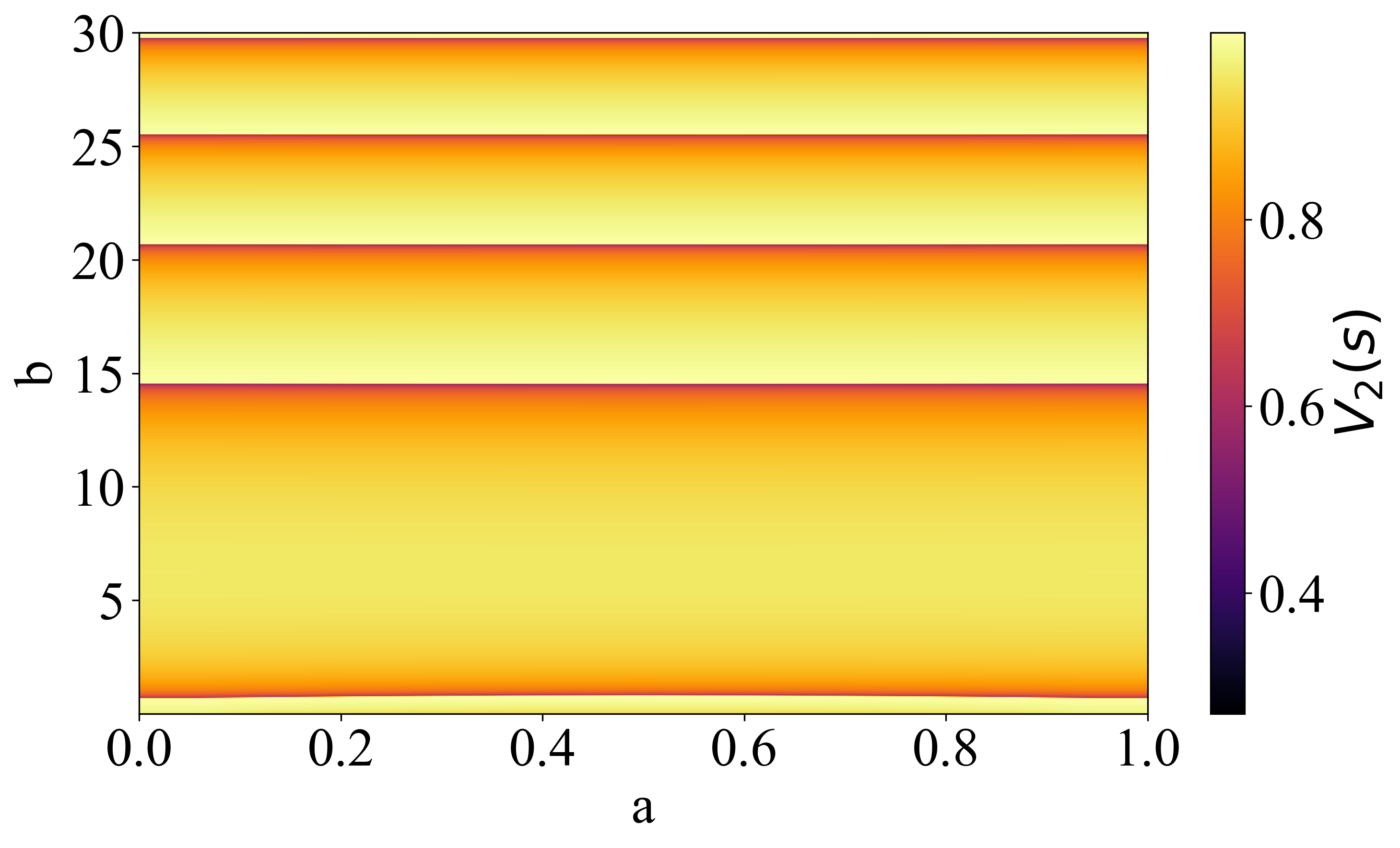}
  \caption{$V_2(s)$ \eqref{V1-V2} as a function of $a$ and $b$.}
  \label{V_2(s)}
\end{center}
\end{figure}

Numerical exploration of the complex plane shows that $V_1(s)$ attains
minimal values close to $\Re(s)=1/2$ and $\Im(s) = 2 \pi$ (Figure \ref{V_1(s)}), while $V_2(s)$ vanishes at imaginary parts corresponding to the known non-trivial zeros of $\zeta(s)$ (Figure \ref{V_2(s)}).
The combined condition $V_1(s)\approx0$ and $V_2(s)\approx0$ therefore provides a convenient numerical indicator for visualizing the location of non-trivial zeros.

We emphasize that this criterion is purely heuristic and visual in
nature.
It does not constitute a proof of the Riemann Hypothesis, nor does it
yield new analytic constraints on the zeros.
Rather, it offers a geometric perspective on how the functional symmetry
of $\zeta(s)$ manifests itself in the neighbourhood of the critical line,
complementing classical analytic criteria such as those of Hardy and Li
\cite{Li1997}.

\subsection{Asymptotic scaling of non-trivial zeros}

The asymptotic distribution of the non-trivial zeros of the Riemann zeta
function is governed by the classical Riemann-von Mangoldt formula,
which states that the number $N(b)$ of zeros with imaginary part in
$(0,b)$ satisfies
\begin{equation}
N(b)
\;\sim\;
\frac{b}{2\pi}
\ln\!\left(\frac{b}{2\pi e}\right),
\label{RVM}
\end{equation}
as $b\to\infty$ \cite{Titchmarsh}.

\begin{figure}[ht]
\begin{center}
  \includegraphics[width=0.5\columnwidth]{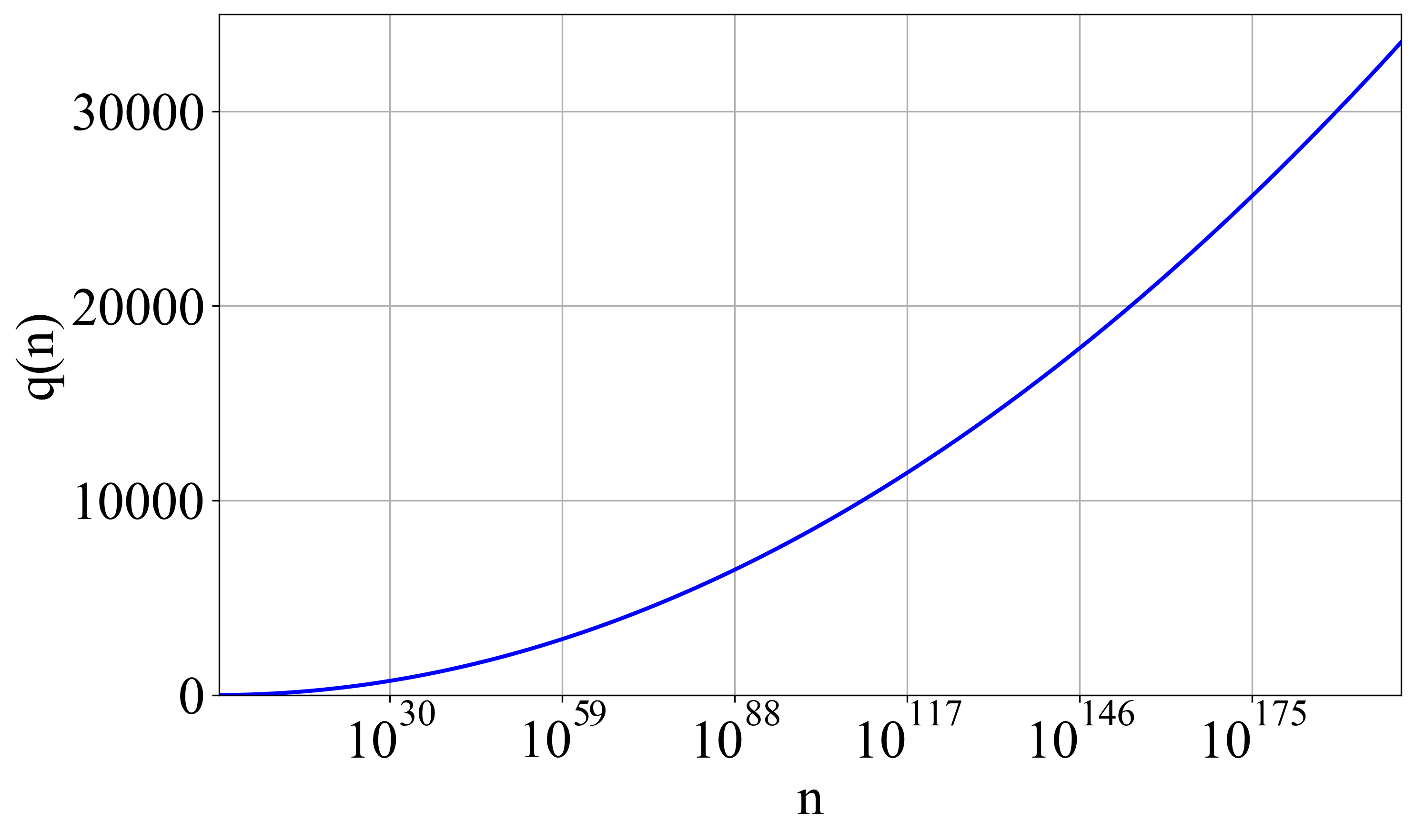}
  \caption{The ratio $q(n)$ \eqref{q(n)_eq}, illustrating
  the relative growth of non-trivial zeros and primes.}
  \label{q(n)}
\end{center}
\end{figure}

Motivated by the comparison between prime statistics and zero statistics,
we introduce the ratio
\begin{equation}
q(n)
\;=\;
\frac{N(n)}{\pi(n)}
\;\sim\;
\frac{\ln(n/2\pi e)\,\ln n}{2\pi},
\label{q(n)_eq}
\end{equation}
which compares the growth rate of non-trivial zeros up to height $n$
with the density of primes below $n$.
This quantity provides a convenient measure of the relative asymptotic scales
of the two distributions (Figure \ref{q(n)}).

Formally inverting the asymptotic relation \eqref{q(n)_eq} yields
\begin{equation}
n(q)
\;\sim\;
\sqrt{2\pi e}\,
\exp\!\left(
\frac{1}{2}
\sqrt{\bigl(1+\ln(2\pi)\bigr)^2+8\pi q}
\right),
\label{n(q)_eq}
\end{equation}

This inversion should be interpreted only at the level of leading-order
asymptotics.
In particular, $n(q=0)=2\pi e\approx17$ provides a rough scale compatible with the appearance of the first
non-trivial zero.
$n(q = 1) \sim 74$ marks the scale at which the expected counts of primes
and zeros become comparable.

Further insight can be obtained by relating the zero-counting function
to the asymptotic growth of the Bernoulli numbers.
Using the approximation
\begin{equation}
|B_{2n}|
\;\sim\;
\sqrt{16\pi n}
\left(\frac{n}{\pi e}\right)^{2n},
\qquad n\gg1,
\end{equation}
one finds the asymptotic representation
\begin{equation}
N(2n)
\;\sim\;
\frac{n}{\pi}
\ln\!\left[
\left(
\frac{|B_{2n}|}{\sqrt{16\pi n}}
\right)^{\!1/2n}
\right],
\label{new_N(n)}
\end{equation}
which is consistent with the Riemann-von Mangoldt formula.
Conversely,
\begin{equation}
|B_{2n}|
\;\sim\;
\sqrt{16\pi n}\,
e^{2\pi N(2n)},
\label{new_B_2n}
\end{equation}
providing an alternative asymptotic description of the growth of the
Bernoulli numbers.

The absence of non-trivial zeros for $b\lesssim14$ is numerically consistent with the fact that $|B_{2n}|<1$ for small $n$, without implying any causal or structural connection.
Using the classical identity
\eqref{zeta-bernoulli}, one further obtains, for $n\gg1$,
\begin{equation}
\zeta(2n)
\;\sim\;
\sqrt{4\pi n}\,
\frac{(2\pi)^{2n}}{(2n)!}\,
e^{2\pi N(2n)},
\label{zeta_N}
\end{equation}
which relates the behaviour of $\zeta(s)$ on the real axis to the global
distribution of its non-trivial zeros.

More refined versions of the zero-counting formula,
including lower-order corrections,
are well known in the literature \cite{Backlund1914,Edwards1974}:
\begin{equation}
N(b) = \frac{b}{2 \pi} \ln \bigg( \frac{b}{2 \pi e} \bigg)
+ \frac{7}{8}
+ S(b)
+ O(b^{-1}),
\label{RVM_2}
\end{equation}
where $S(b) = \arg \zeta(1/2+ib) / \pi$.

In particular, they allow one to approximate the imaginary part $b_k$
of the $k$-th non-trivial zero via the Lambert $W$-function, for $k \gg 1$:
\begin{equation}
b_k
\;\sim\;
\frac{2\pi\,(k-11/8)}{W_0\!\left[(k-11/8)/e\right]}.
\end{equation}

Substituting it into the truncated Hadamard product
\eqref{zeta(s)_new}, we obtain a truncated explicit representation of $\zeta(s)$,
\begin{equation}
\zeta(s)
\;\sim\;
f(s)
\prod_{k=1}^{n_{\max}}
g(s,k)\,h(k)^s
\frac{1-s/s_k}{e^{-s(1/2k+1/s_k)}},
\label{zeta(s)_new_2}
\end{equation}
where $\Im(s) \gg 1$, $s_k=1/2+ib_k$ and
$n_{\max}=\lfloor (b/2\pi)\ln(b/2\pi e)\rfloor+1$.
This expression emphasizes how the contribution of individual zeros
builds up the global structure of $\zeta(s)$.

\subsection{Oscillatory structure of the Euler product}

We analyze the Euler product representation of the Riemann zeta function,
\begin{equation}
\zeta(s) = \prod_{p} \frac{1}{1 - p^{-s}},
\end{equation}
by studying the local behavior of each Euler factor in the complex plane.
Let $s = a + ib$, and define
\[
z_p = p^{-a}, \qquad \theta_p = b \ln p.
\]
The modulus and phase of $(1 - p^{-s})^{-1}$ are given explicitly by:
\begin{equation}
r_p(s) = \left( 1 - 2 z_p \cos\theta_p + z_p^2 \right)^{-1/2}, \qquad
\varphi_p(s) = \arctan^2 \left(
\frac{z_p \sin\theta_p}{1 - z_p \cos\theta_p}
\right)
\label{r_p_def}
\end{equation}
where $\varphi_p(s) \in (-\pi,\pi)$.

Since $\cos\theta_p \in [-1,1]$, the modulus $r_p(s)$ satisfies the bounds
\begin{equation}
(r_p)_{min} = \frac{1}{1 + z_p} \le r_p(s) \le \frac{1}{1 - z_p} = (r_p)_{max},
\label{r_p_bounds}
\end{equation}

Formally averaging over a uniform distribution of phases $\theta_p$ yields
\begin{equation}
\langle r_p \rangle \approx \frac{1}{1 - z_p^2},
\label{r_p_average}
\end{equation}
which is strictly larger than $1$ for $a > 0$.
This observation highlights that the convergence properties of the Euler product
cannot be inferred from amplitude considerations alone, but crucially depend on
the collective interaction of the oscillatory phases.

\begin{figure}[ht]
\begin{center}
  \includegraphics[width=0.5\columnwidth]{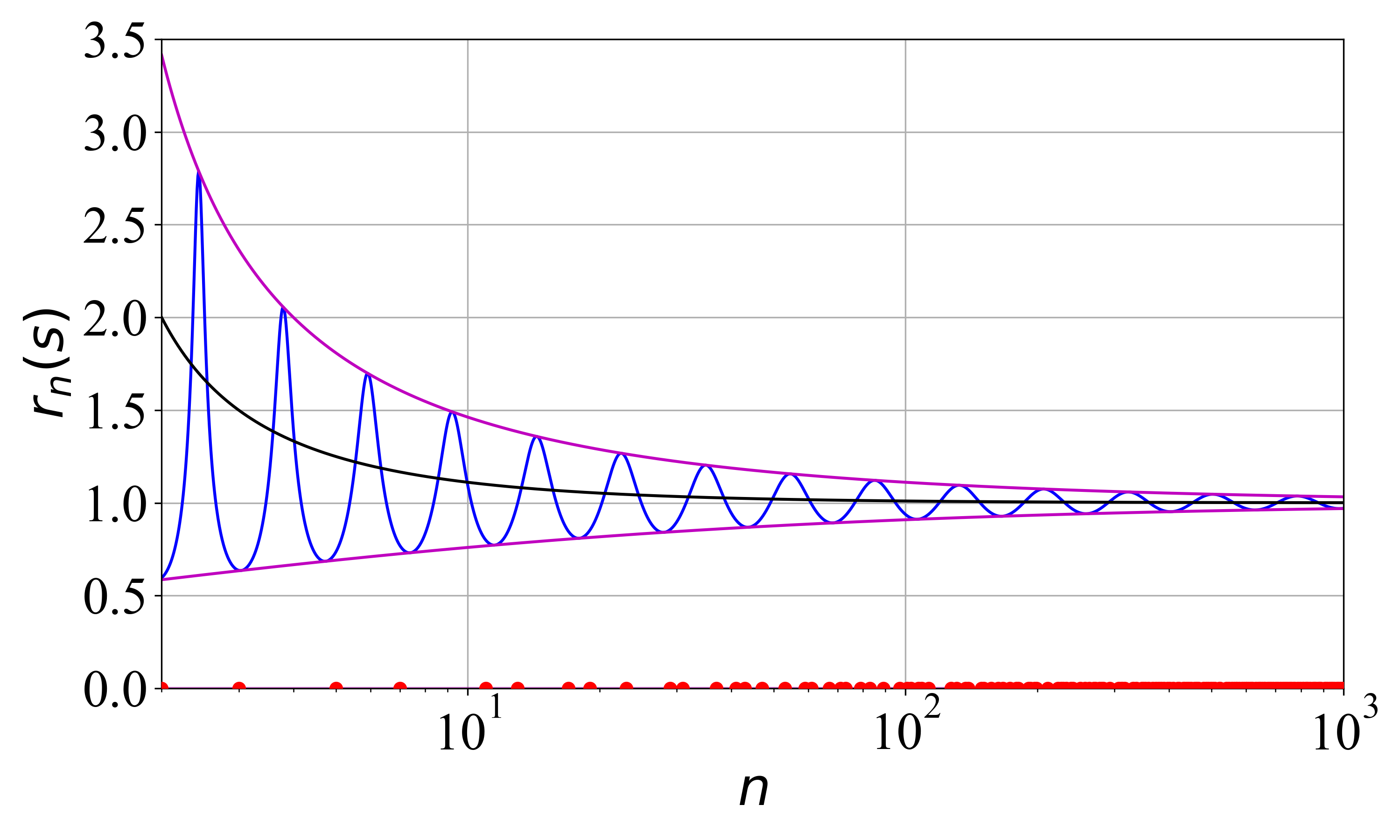}
  \caption{$r_n(0.5 + 14.135  \ i)$ \eqref{r_p_def} (blue), $(r_p)_{min}$ and $(r_p)_{max}$ \eqref{r_p_bounds} (magenta), and $\langle r_p \rangle$ \eqref{r_p_average} (black). Red points indicate the positions of prime numbers.}
  \label{r(n)}
\end{center}
\end{figure}

Figure~\ref{r(n)} illustrates the behavior of $r_n(s)$ evaluated at the first
non-trivial zero $s = 1/2 + 14.135\, i$.
The observed alignment between local minima of $r_n(s)$ and the distribution of prime numbers (indicated by red points) suggests that, at the imaginary parts corresponding to non-trivial zeros,
the Euler factors exhibit a highly structured interference pattern.
This phenomenon is consistent with heuristic models in which cancellations among
local oscillations play a central role in the emergence of zeros of $\zeta(s)$.

\begin{figure}[ht]
\begin{center}
  \includegraphics[width=0.5\columnwidth]{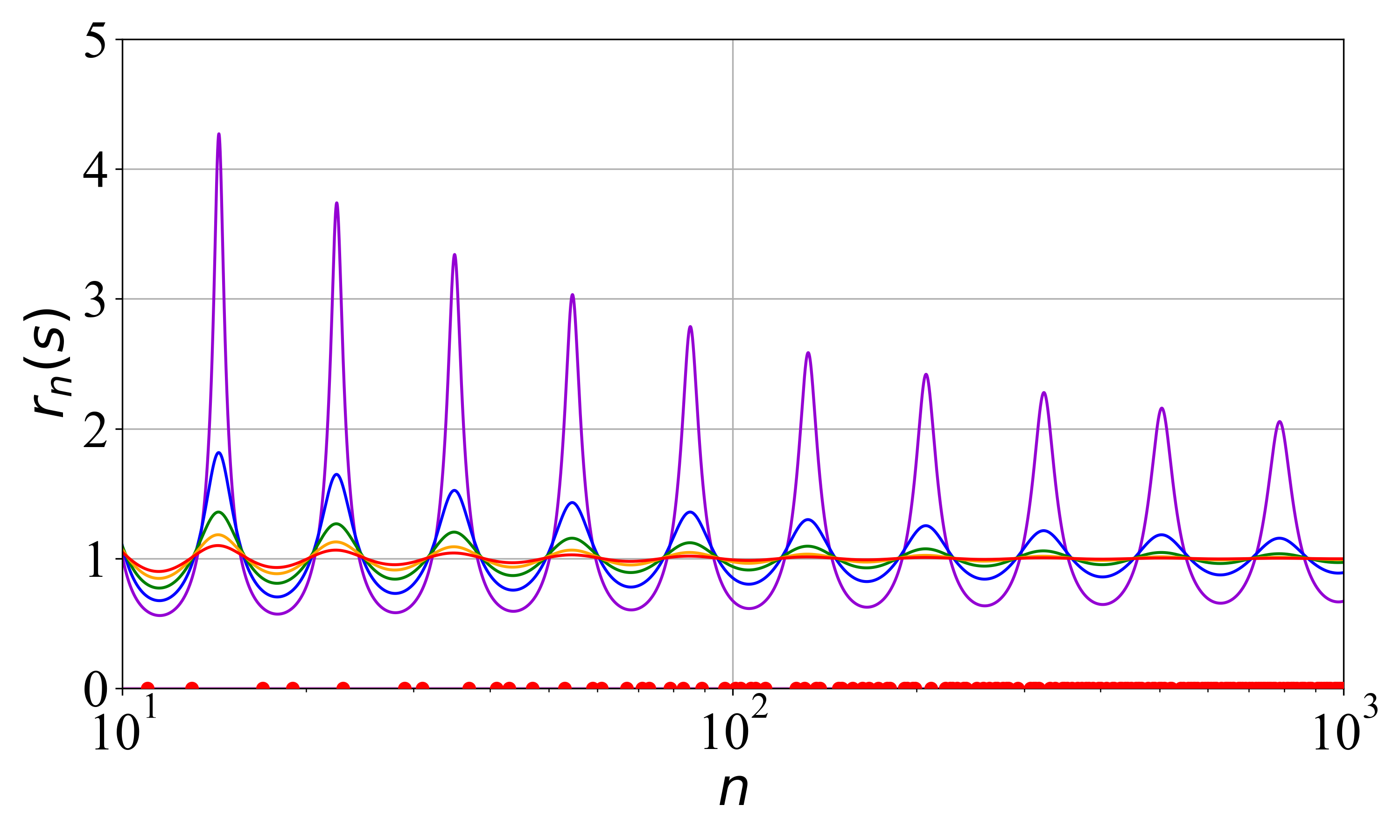}
  \caption{$r_n(s)$ \eqref{r_p_def} for $b = 14.135$: $a = 0.1$ (violet), $a = 0.3$ (blue), $a = 0.5$ (green), $a = 0.7$ (orange), $a = 0.9$ (red). Red points indicate the positions of prime numbers.}
  \label{r(n)_2}
\end{center}
\end{figure}

The oscillatory behavior of $r_p(s)$ is governed by the phase $\theta_p = b \ln p$,
leading to fluctuations whose amplitude is
\begin{equation}
\Delta r_p = r_{p,\max} - r_{p,\min}
= \frac{2 p^{a}}{p^{2a} - 1}.
\label{delta_r}
\end{equation}
In particular, on the critical line $a = 1/2$ and for large primes $p$, one has 
$ \Delta r_p \sim 2/\sqrt{p} $,
indicating a slow decay of oscillation amplitudes with increasing $p$ (Figure \ref{r(n)_2}).

For primes $p \le n$, the number of oscillations of $r_p(s)$ up to scale $n$ 
can be heuristically estimated as
\begin{equation}
X(b,n) \sim \frac{b}{2\pi} \ln\!\left( \frac{n}{2} \right),
\end{equation}
which mirrors the logarithmic scaling appearing in the Riemann-von Mangoldt
formula and is consistent with Montgomery’s analysis of the pair correlation
of non-trivial zeros \cite{Montgomery1973}.
Comparing this scaling with the Riemann-von Mangoldt formula, one obtains
\begin{equation}
N(b) \sim X(b,n) - \frac{b}{2 \pi} \ln \!\Big( \frac{\pi e n}{b} \Big),
\end{equation}
up to lower-order terms, with $N(b)=0$ for $n = 2\pi e$.

It is then natural to ask how many extrema of $r_p(s)$ can effectively
contribute to the truncated Euler product.
If $X(b,n) < K(n)$, the number of primes exceeds the number of oscillations,
while if $X(b,n) > K(n)$, the density of primes limits the number of extrema
that may contribute independently to $\zeta(s)$.

Equating $X(b,n) = K(n)$ suggests a characteristic crossover scale for the imaginary part,
\begin{equation}
b_{\mathrm{lim}}(n) = \pm \, \frac{2\pi K(n)}{\ln(n/2)}
\sim \pm \, \frac{2\pi n}{\ln(n/2)\ln n}.
\label{b_lim_eq}
\end{equation}

The function $b_{\mathrm{lim}}(n)$ attains its minimum value $\sim 2\pi e$
for $n \approx 9-13$, a scale numerically comparable 
to the height of the first non-trivial zero, consistently with Equation \eqref{q(n)_eq}.

In particular, for the first non-trivial zero, one has $b < b_{\mathrm{lim}}(n)$
for all $n$, indicating that no oscillatory extrema are a priori excluded
from contributing in the early stages of the Euler product.

We consider the average modulus of the Euler factors
over the first $K(k)$ primes,
\begin{equation}
M_k(s) = K(k)^{-1} \sum_{p \le k} r_p(s),
\label{M_b}
\end{equation}
where $K(k)$ denotes the number of primes not exceeding $k$.
For $k=50$, numerical evaluation over the first $29$ non-trivial zeros
($b<100$) yields
\begin{equation}
\min_{p \le k} r_p(b) \sim 0.58-0.71, \qquad
M_k(b) \sim 0.92-0.97, \qquad
\max_{p \le k} r_p(b) \sim 1.17-1.47.
\end{equation}

\begin{figure}[ht]
\begin{center}
  \includegraphics[width=0.5\columnwidth]{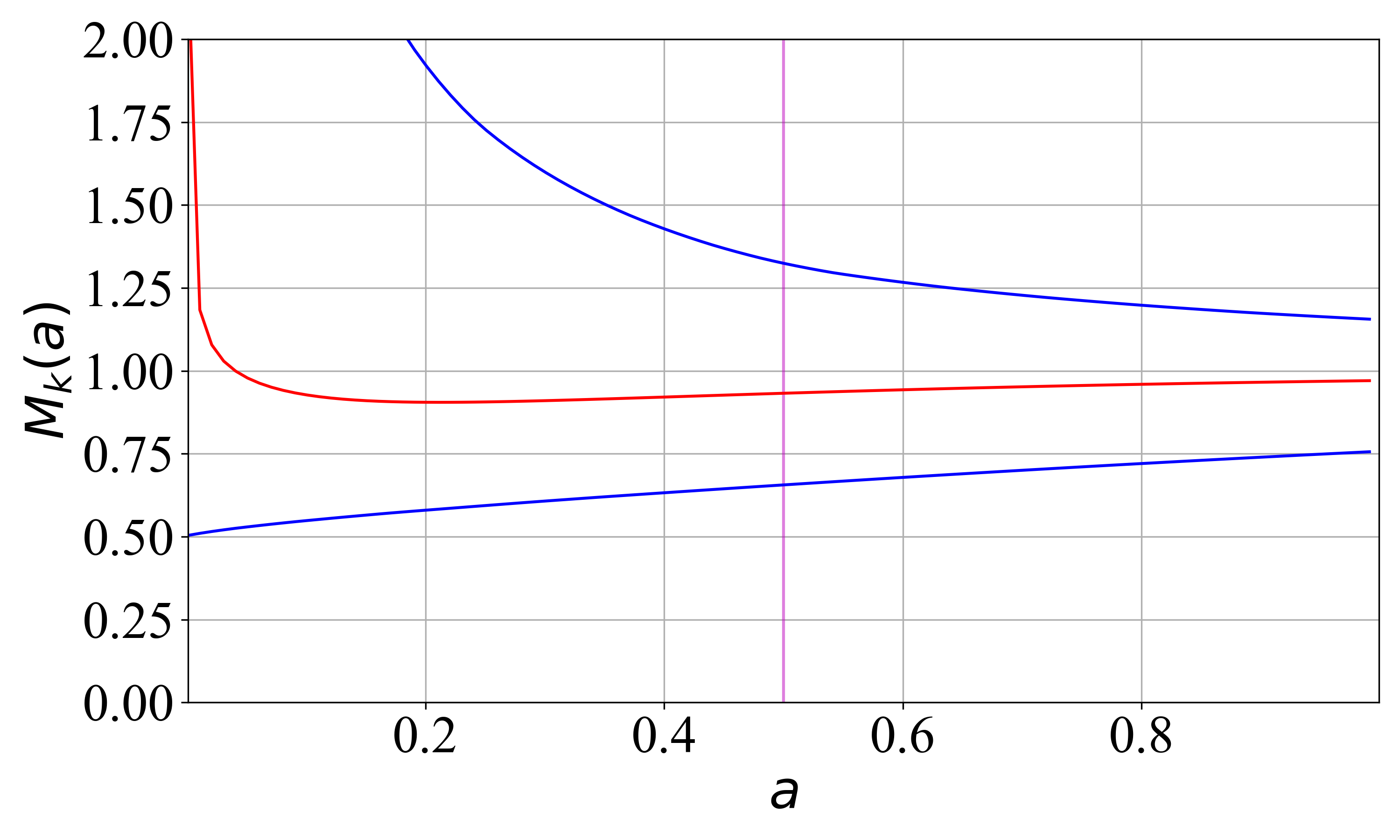}
  \caption{Average values of $M_k(b)$ (red), $\min_{p \leq k} r_p(b)$ and $\max_{p \leq k} r_p(b)$ (blue), for $k = 50$ and $b < 100$. As reference, the critical line ($a = 1/2$) is shown in magenta.}
  \label{M(a)_plot}
\end{center}
\end{figure}

Figure~\ref{M(a)_plot} displays the averages of $M_k(b)$ together with the
corresponding minima and maxima as functions of $a$, for fixed $k=50$ and
$b<100$.
On the critical line $a=1/2$, the theoretical bounds
$(1+2^{-1/2})^{-1} \approx 0.586 $ and
$(1-2^{-1/2})^{-1} \approx 3.413 $
are consistent with the observed numerical range.
These results suggest that, even near non-trivial zeros, the local Euler
factors remain bounded and moderately distributed, reinforcing the widely held view that zeros of $\zeta(s)$ arise from delicate global cancellations among Euler factors \cite{Soundararajan2009}.

\subsection{Fluctuations and stability of the Euler product}

The dependence of the local Euler factors $r_p(s)$ on the prime variable $p$
is governed by
\begin{equation}
\frac{\partial r_p}{\partial p}
=
- \frac{z_p}{p} \, r_p(s)^{-3}
\bigl[
a(\cos\theta_p - z_p) + b \sin\theta_p
\bigr],
\end{equation}
where $z_p = p^{-a}$ and $\theta_p = b \ln p$.
Here $p$ is treated as a continuous variable, as customary in heuristic
analyses of Euler products.

For any $a>0$, one has $r_p(s)\to 1$ as $p\to\infty$,
while both the amplitude and the frequency of the oscillatory component decrease.
As a consequence, the initial factors of the Euler product dominate the global
behavior of $|\zeta(s)|$, consistently with previous analyses of truncated Euler
products \cite{GonekHughesKeating2007}.

Using the expansion
\begin{equation}
(1-u)^{-1/2}
=
1 + \frac{1}{2}u + \frac{3}{8}u^2 + O(u^3),
\end{equation}
one obtains, for $p\gg1$ and $a>0$,
\begin{equation}
r_p(s)
\sim
1 + z_p \cos\theta_p
+ z_p^2\!\left(\frac14 + \frac34 \cos(2\theta_p)\right).
\end{equation}
This decomposition naturally separates a not-oscillatory contribution $1 + z_p^2/4$
and two harmonic components:
\begin{equation}
H_{1,p} = z_p \cos\theta_p,
\qquad
H_{2,p} = \frac34 z_p^2 \cos(2\theta_p).
\end{equation}

To assess convergence properties, it is therefore natural to study
the cumulative effect of these oscillatory terms on
$\sum_{p\le n} \log r_p(s)$.

Formally differentiating $\theta_p$ with respect to $p$ yields the local angular frequency
$\omega(p)=|b|/p$, implying an oscillation period
$\Delta p\sim 2\pi p/|b|$.
According to the Prime Number Theorem, the number of primes within such a period is
asymptotically
\begin{equation}
N(p)\sim 2\pi p/(|b|\ln p).
\label{N(p)_period}
\end{equation}

During a period of oscillation, the main contribution of $H_{1,p}$ and $H_{2,p}$ to $|\zeta(s)|$ is $\sum_{j=1}^{N(p)} z_j \cos \theta_j$.
Assuming the phases $\theta_p$ to be uniformly distributed modulo $2\pi$,
the variance of $\cos\theta_p$ is $1/2$.
Therefore, the standard deviation of the cumulative harmonic contribution is
\begin{equation}
\sigma_p(s) \sim \sqrt{\frac{1}{2} \sum_{j=1}^N z_j^2}
\end{equation}

Assuming quasi-independence of phases over a period and approximating
$z_j \sim p^{-a}$ within the interval, we obtain
\begin{equation}
\sigma_p(s)
\sim
\sqrt{\frac{\pi}{|b|\ln p}}\; p^{1/2-a}.
\label{sigma_p}
\end{equation}

In general, every period of oscillation of $r_p(s)$ can heuristically 
induce multiplicative fluctuations of order $e^{\pm \sigma}$.
On average, every prime contributes to the Euler product with a factor $\exp(\pm p^{-a}/\sqrt{2})$.

For $a>1/2$, $\sigma_p(s)\to0$ as $p\to\infty$, so that fluctuations of the Euler
product become negligible and $|\zeta(s)|$ stabilizes.
For $a<1/2$, $\sigma_p(s) \to \infty$ as $p \to \infty$, indicating global instability of the product.

The Euler product presents stable fluctuations if $\ln \sigma_p(s) < 0$ for $p \to \infty$.
In particular, the fastest convergence of $\zeta(s)$ to $0$ happens for the highest possible $\ln \sigma_p(s)$ lower than $0$, suggesting the following condition for $a$:
\begin{equation}
a = \frac{1}{2} - \frac{1}{2 \ln p} \ln \bigg( \frac{|b|}{\pi} \ln p \bigg)
\label{a_constraint}
\end{equation}
which, for $p \to \infty$, becomes $a = 1/2$.

\begin{figure}[ht]
\begin{center}
  \includegraphics[width=0.5\columnwidth]{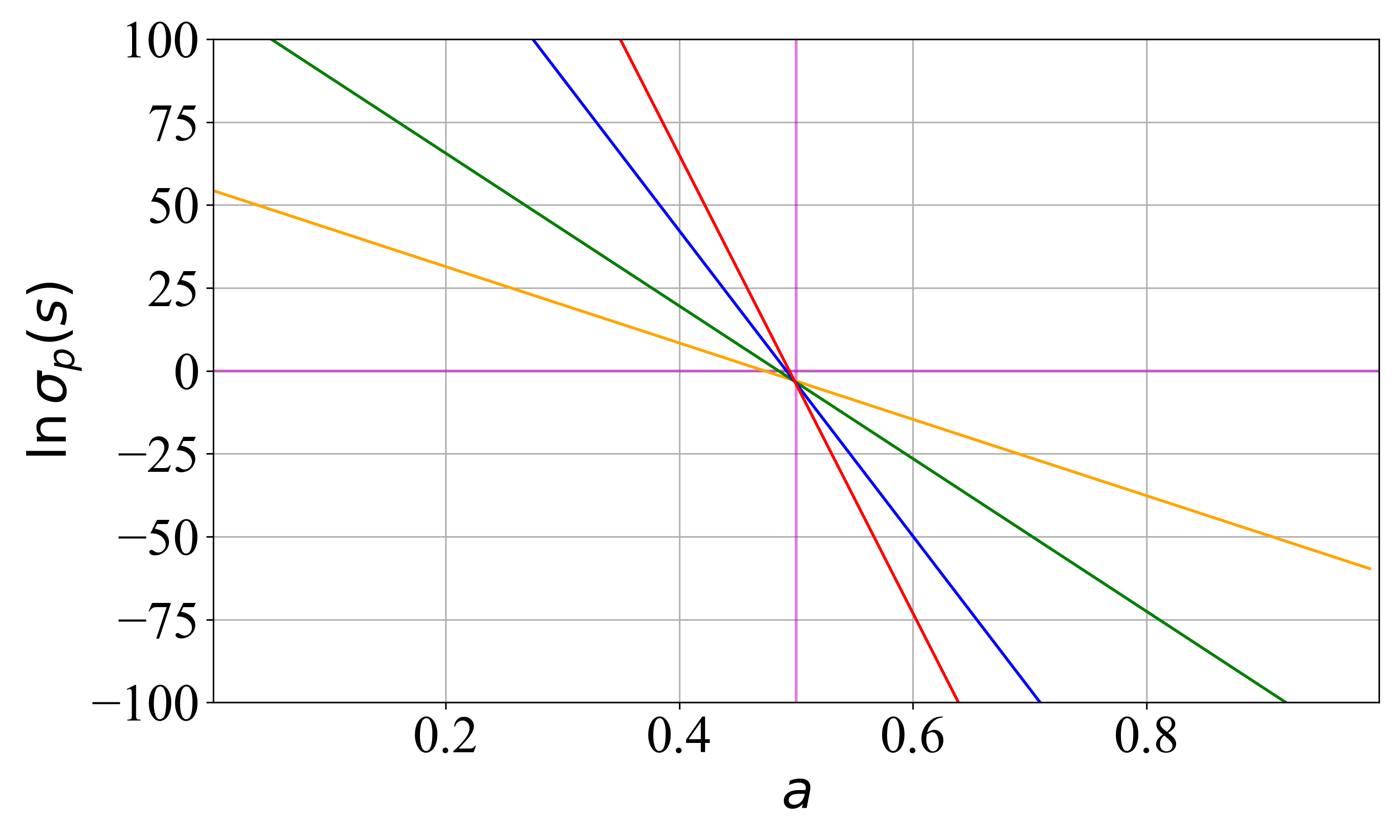}
  \caption{$\ln \sigma_p(s)$ \eqref{sigma_p} for $b = 14.135$: $p = 10^{50}$ (orange), $p = 10^{100}$ (green), $p = 10^{200}$ (blue), $p = 10^{300}$ (red). As reference, $a = 1/2$ and $\ln \sigma_p(s) = 0$ are shown in magenta.}
  \label{ln_sigma_p(s)}
\end{center}
\end{figure}

The critical line therefore corresponds to the slowest possible decay of
$\sigma_p(s)$ compatible with stability,
\begin{equation}
\sigma_p(1/2+b)\sim \sqrt{\frac{\pi}{|b|\ln p}},
\end{equation}
reflecting a precise balance between oscillation amplitude and prime density (Figure \ref{ln_sigma_p(s)}).

While this argument is heuristic and does not constitute a proof of the
Riemann Hypothesis, it provides a structural heuristic for why the critical line plays a
distinguished role in the fluctuation balance of the Euler product,
in agreement with classical insights \cite{Selberg1946} and modern analyses of zeta-function
fluctuations \cite{Harper2014}.

\subsection{Local minima of truncated Euler products}

From the definition \eqref{r_p_def},
\begin{equation}
r_p(s) \leq 1 \ \Leftrightarrow \ 2\cos\theta_p \le \frac{1}{2}p^{-a},
\end{equation}
which is satisfied for $b$ belonging to the union of intervals
\begin{equation}
b \in \bigcup_{k\in\mathbb{Z}}
\left[
\frac{2\pi k+\alpha}{\ln p},
\frac{2\pi(k+1)-\alpha}{\ln p}
\right],
\label{b_interval}
\end{equation}
where $\alpha=\arccos\!\left(p^{-a}/2 \right)$.

These intervals correspond to phases for which the local Euler factor
induces a contraction rather than an amplification of the truncated product.
These intervals have spacing $\Delta b = 2\pi/\ln p$ and width
$2(\pi-\alpha)/\ln p$, so that the fraction of admissible values is
$f_p(a)=1-\alpha/\pi$.
In particular, for $a=1/2$ one finds $f_p(a)\to 1/2$ as $p\to\infty$,
reflecting an asymptotic balance between amplifying and damping phases.

When a value of $b$ satisfies the condition $r_p(s)<1$ for a substantial
proportion of primes, the truncated Euler product is more likely to exhibit
a local suppression of its modulus.
This motivates the introduction of the averaged truncated modulus
$M_k(b)$ defined in Equation \eqref{M_b}.
Numerically, the imaginary parts of non-trivial zeros are observed to lie
systematically near local minima of $M_k(b)$, consistently with previous studies of truncated
Euler products \cite{GonekHughesKeating2007}.

Treating $b$ as a continuous parameter and keeping the truncation index $k$
fixed, differentiation yields
\begin{equation}
\frac{\partial M_k}{\partial b} = K(k)^{-1} \sum_{p \leq k} \frac{z_p \ \theta_p}{b} \frac{\sin \theta_p}{D_p^{3/2}}, \qquad
\frac{\partial ^2M_k}{\partial b^2} = K(k)^{-1} \sum_{p \leq k} \frac{z_p \ \theta_p^2}{b^2} \bigg( \frac{\cos \theta_p}{D_p^{3/2}} - \frac{3 z_p \sin^2 \theta_p}{D_p^{5/2}} \bigg),
\label{M_k_der}
\end{equation}
where $D_p = 1 - 2z_p\cos\theta_p + z_p^2$.

From these expressions one is led to the following heuristic numerical
criteria for the localization of minima:
\begin{subequations}
\begin{align}
&K(|b|)^{-1}\sum_{p\le |b|} r_p(s) < 1 + |b|^{-1},
\\
&K(|b|)^{-1}\Bigl|\sum_{p\le |b|} \sin\theta_p\Bigr| < |b|^{-1},
\\
&K(|b|)^{-1}\sum_{p\le |b|} \cos\theta_p < |b|^{-1},
\end{align}
\label{b_constraints}
\end{subequations}
where $K(|b|)$ denotes the number of primes not exceeding $|b|$ 
and $|b|^{-1}$ represents a numerical tolerance.

\begin{figure}[ht]
\begin{center}
  \includegraphics[width=0.5\columnwidth]{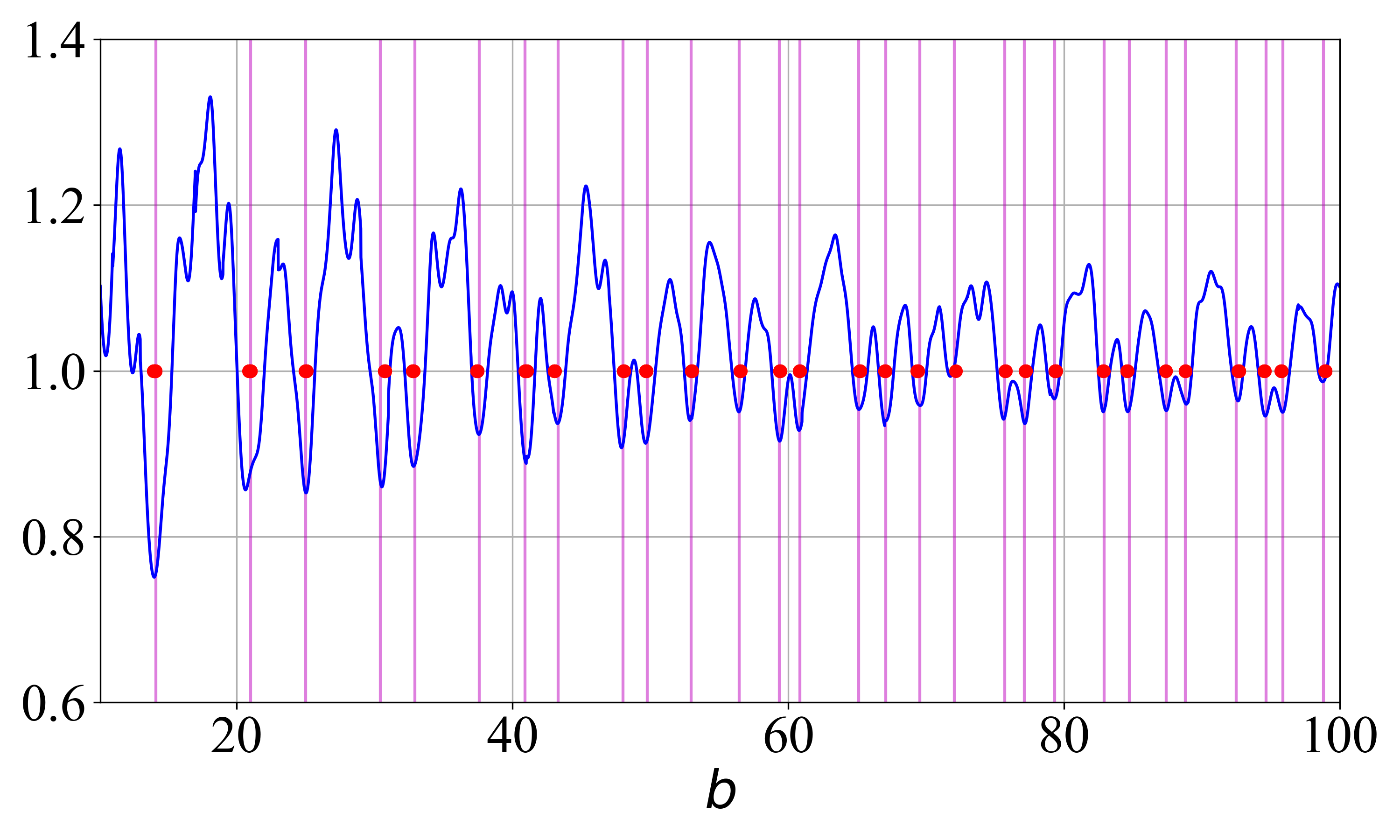}
  \caption{$M_k(b)$ \eqref{M_b} (blue) and the imaginary parts of the non-trivial zeros of $\zeta(s)$ (magenta), for $k = b$, $a = 1/2$ and $b < 100$. Red dots indicate the values of $b$ for which the 3 constraints \eqref{b_constraints} are satisfied simultaneously.}
  \label{M(b)_plot}
\end{center}
\end{figure}

Empirically, values of $b$ satisfying these conditions 
provide a reliable numerical proxy for the imaginary parts of non-trivial zeros
using only the first $K(|b|)$ primes and without explicit evaluation of the phase of $\zeta(s)$
(Figure~\ref{M(b)_plot}).

To quantify the prevalence of local damping, we define $F_k(s)$ as the proportion of primes for which the local Euler factor has modulus smaller than unity:
\begin{equation}
F_k(s) = \frac{1}{K(k)} \sum_{p\le k} \mathbf{1}_{\{r_p(s)<1\}},
\label{F_k(s)}
\end{equation}

\begin{figure}[ht]
\begin{center}
  \includegraphics[width=0.5\columnwidth]{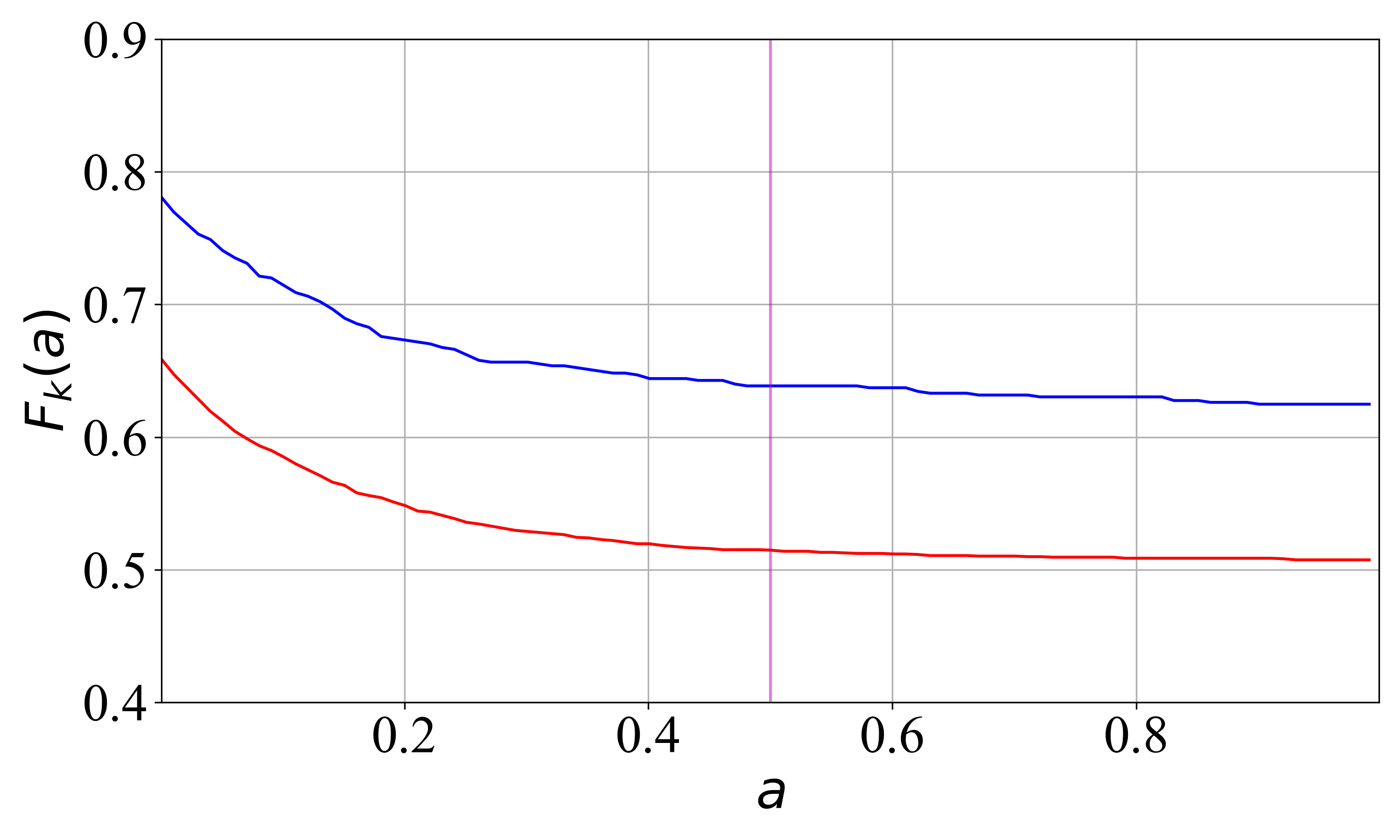}
  \caption{$\langle F_k(a) \rangle$ \eqref{F_k(s)} for $b < 100$ (red) and for the first 29 non-trivial zeros (blue). As reference, the critical line ($a = 1/2$) is shown in magenta.}
  \label{F_k(a)_plot}
\end{center}
\end{figure}

Numerical experiments show that $F_k(s)$ is systematically larger near
local minima of $M_k(b)$ and when $a < 1/2$.
In particular, for $k=100$, $F_k(s) \sim 62.5 - 78.2 \%$ for values of $b$ corresponding to non-trivial zeros and $F_k(s)$ is $\sim 50.3 - 66.7 \%$ for generic values of $b$ (Figure \ref{F_k(a)_plot}).
These values are stable under moderate variations of $k$.

In summary, the imaginary parts of non-trivial zeros are characterized by a
concurrent minimization of $M_k(b)$ and an enhanced proportion of primes
satisfying $r_p(s)<1$.
As $p\to\infty$, the fastest stable suppression of the Euler product
is asymptotically compatible only with $a=1/2$ within this heuristic framework, 
in agreement with the critical-line phenomenon predicted by the Riemann Hypothesis.

\subsection{Pair correlation and Euler-product oscillations}

The imaginary-component criterion introduced above admits a natural
interpretation in terms of the pair correlation of the non-trivial zeros of
$\zeta(s)$.
Montgomery's pair correlation conjecture \cite{Montgomery1973}, supported by
extensive numerical evidence \cite{Odlyzko1987}, shows that the imaginary
parts of the zeros exhibit strong local correlations and short-range
repulsion, governed by an underlying oscillatory structure.

In the present framework, an analogous oscillatory behavior emerges through
the prime-dependent phases $\theta_p = b \ln p$ appearing in the Euler product.
Although these phases are indexed by primes rather than zeros, their collective
interference governs the fine structure of truncated Euler products in the
imaginary direction.

Local minima of the averaged modulus $M_k(b)$ 
arise when the oscillations of $r_p(s)$ across different primes exhibit a
statistically coherent alignment.
These values of $b$ correspond to regions where $|\zeta(1/2+ib)|$ is
statistically likely to be small.

From this perspective, the minima of $M_k(b)$ do not identify individual zeros
directly, but rather mark statistically preferred locations determined by the
same interference mechanisms that underlie zero repulsion.

While Montgomery's theory probes correlations between the zeros themselves,
the present criterion probes correlations between prime-indexed oscillatory
contributions in the Euler product.
Although these objects are different, both descriptions capture complementary
manifestations of a common harmonic structure linking primes and zeros.

This connection is consistent with semiclassical and random-matrix approaches
in which Euler products are interpreted as generating functions for spectral
statistics \cite{BogomolnyKeating1996},
and provides a natural statistical explanation for why truncated Euler products
are numerically effective in localizing the imaginary parts of non-trivial zeros.

\newpage

\section*{Conclusions}

Unlike previous approaches, introducing density-based structural frameworks for prime numbers and prime gaps, we:
\begin{itemize}
\item{explicitely unified primality, coprimality and prime pairs properties};
\item{showed a structural tension between HL, Cramér, and PNT predictions, and obtained quantitative estimates on the rarity of extreme gaps};
\item{derived non-conjectural bounds for Goldbach-type additive representations};
\item{analyzed the oscillatory structure of truncated Euler products, providing a coherent interpretation of the critical line and the imaginary part of non-trivial zeros}.
\end{itemize}

Despite this, several fundamental questions
remain open. In particular, a fully rigorous derivation of extreme gap
distributions, a comprehensive statistical
theory of additive representations within this normalized setting, 
and a precise analytic formulation of Euler product stability
sufficient to address the Riemann Hypothesis remain to be developed. 
Addressing these problems may then provide a unified structural basis for
probabilistic, combinatorial, and analytic approaches to some of the central
questions of number theory.

\bibliographystyle{amsplain}
\bibliography{bibliography}

\end{document}